%% file: main.tex
\documentclass[11pt,letterpaper]{article}
\usepackage[T1]{fontenc}
\usepackage[utf8]{inputenc}
\usepackage{lmodern}
\usepackage[letterpaper,margin=1in]{geometry}
\usepackage{amsmath,amssymb,amsthm,mathtools}
\usepackage{microtype,enumitem,xcolor}
\usepackage{hyperref}
\hypersetup{colorlinks=true,filecolor=blue!45!black,linkcolor=blue!45!black,citecolor=blue!45!black,urlcolor=blue!45!black,pdftitle={Integral saturation of adjacency-degree lattices},pdfauthor={Zhipeng Lu}}
\setlist{itemsep=3pt,topsep=5pt}
\newtheorem{theorem}{Theorem}[section]
\newtheorem{lemma}[theorem]{Lemma}
\newtheorem{proposition}[theorem]{Proposition}
\newtheorem{corollary}[theorem]{Corollary}
\newtheorem*{theorem*}{Theorem}
\newtheorem*{lemma*}{Lemma}
\newtheorem*{proposition*}{Proposition}
\theoremstyle{definition}
\newtheorem{definition}[theorem]{Definition}
\newtheorem{example}[theorem]{Example}

\theoremstyle{remark}
\newtheorem{remark}[theorem]{Remark}

\numberwithin{equation}{section}
\allowdisplaybreaks[1]
\title{Integral saturation of adjacency--degree lattices}
\author{Zhipeng Lu\\
\normalsize Shenzhen MSU--BIT University\\
\normalsize Guangdong Laboratory of Machine Perception and Intelligent Computing\\
\normalsize Shenzhen 518172, Guangdong, China\\
\normalsize\texttt{zhipeng.lu@hotmail.com}}
\date{}
\begin{document}
\hypersetup{pageanchor=false}
\begin{titlepage}
\maketitle
\thispagestyle{empty}
\begin{abstract}\normalsize
Let \(A\) be the adjacency matrix of \(G\sim G(n,1/2)\), let
\(D=\operatorname{diag}(A\mathbf1)\), and let
\(\Lambda_G=\mathbb Z\langle A,D\rangle\mathbf1\).
We prove
\[
 \Pr[\Lambda_G\ne\mathbb Z^n]\le n^{-3/20+o(1)}.
\]
Thus two dependent graph operators generate the entire vertex lattice
from the constant vector, simultaneously over every prime field.
Adjacency alone has a modulo-two rank obstruction for every simple
graph of order greater than one, even though its rational word span
is typically full.

The proof first constructs words containing at most two adjacency
operations whose integer span contains \(\Delta_G^3\mathbb Z^n\),
where \(\Delta_G\) is the product of the pairwise differences
between distinct degree values. This confines all possible arithmetic
obstructions to a finite prime range on a single event. We then
remove the remaining obstructions using conditional alternating
forms, symmetric random-matrix code estimates, and uniform counts
of degree-graded invariant subspaces.

We connect the integer lattice to graph reconstruction through
pointed moments. Saturation forces every pointed orthogonal AD
intertwiner to be a vertex permutation. The short-word
construction also shows that a fixed polynomial-size family of
ordinary caterpillars, of diameter at most seven and order
\(O(\sqrt{n\log n})\), identifies almost every graph among all
finite simple graphs. For nonsaturated full-rank lattices, we give
an exact description of the moment fibre by admissible unimodular
overlattices; at index two, the only possible ambiguity is a
specified four-vertex Godsil--McKay switch.
\end{abstract}
\end{titlepage}
\setcounter{page}{1}
\hypersetup{pageanchor=true}
\input{introduction}
\input{integer_construction}
\input{prime_mechanisms}
\input{reconstruction}

\label{paper:last-main-page}
\clearpage
\phantomsection\label{paper:submission-proof-index}
\paragraph{Reading and proof guide.}
Sections~\ref{sec:introduction}--\ref{paper:reconstruction} present the
main results and proof framework in the first
\pageref{paper:last-main-page} numbered pages.
Appendices~\ref{sec:foundations}--\ref{sec:caterpillar-reconstruction}
give the finite framework, the full integer construction and final
saturation argument, and the caterpillar proof.
The longer local-prime proofs and the remaining arithmetic theory
are collected in the supplementary proofs beginning on
page~\pageref{paper:additional-proofs}.
All proofs and references are contained in this PDF.
\appendix
\input{article_foundations}
\input{article_random_large_primes}
\input{article_random_integral}
\input{article_caterpillar_reconstruction}

\clearpage
\phantomsection\label{paper:additional-proofs}
\section*{Supplementary proofs}
\addcontentsline{toc}{section}{Supplementary proofs}
The main exposition and Appendices A--C can be read independently of
the order of the following technical proofs.
Sections~\ref{modtwo:section}--\ref{sec:random-growing-primes} prove
the local estimates used to complete Theorem~\ref{integral:random}.
Sections~\ref{sec:arith-intertwiners}--\ref{sec:arith-nonsaturated}
give the full moment-fibre theory and the index-two classification.
Sections~\ref{sec:implications}--\ref{app:integer-certificates}
contain further directions, exact computation, and the integer certificates.
\setcounter{section}{0}
\renewcommand{\thesection}{S\arabic{section}}
\renewcommand{\theHsection}{supplement.\arabic{section}}
\input{article_random_mod2}
\input{article_random_fixed_odd}
\input{article_random_growing_primes}
\input{article_arithmetic_core}
\input{article_implications}
\input{article_computation_appendix}
\input{article_certificates_appendix}
\clearpage
\bibliographystyle{alpha}
\bibliography{references}
\paragraph{AI methodology.}
ChatGPT and OpenAI Codex assisted with the literature review, the
organization and drafting of the exposition, and checks of mathematical
arguments and supporting verification code. The author takes full
responsibility for the statements, proofs, and references. The finite
computational checks support the exposition; they are not formal
verification of the probabilistic arguments.
\end{document}

%% file: introduction.tex
\section{Introduction}\label{sec:introduction}
Starting from the constant function on a graph, one can repeatedly sum
over neighbours or multiply by the vertex degree. Over a field, these
operations produce a cyclic linear representation. Over the integers,
they pose a stronger question: can every individual vertex indicator
be generated without division?

Let \(A\) be the adjacency matrix of a finite simple graph \(G\),
\(e=\mathbf1\), and \(D=\operatorname{diag}(Ae)\). The
\emph{adjacency--degree word lattice} is
\[
 \Lambda_G=\mathbb Z\langle A,D\rangle e
   =\operatorname{span}_{\mathbb Z}\{w(A,D)e:w\text{ a word}\}.
\]
The algebra is unital, and its two generators need not commute.
Equivalently, \(\Lambda_G\) is the smallest subgroup containing \(e\)
and stable under \(A,D\).
We call \(G\) \emph{AD-full} if
\(\mathbb Q\Lambda_G=\mathbb Q^n\), and \emph{integrally AD-saturated}
if \(\Lambda_G=\mathbb Z^n\).

The distinction is substantial even for random graphs. The ordinary
walk matrix
\[ [e,Ae,\ldots,A^{n-1}e] \] is nonsingular with probability
tending to one in \(G(n,1/2)\)~\cite{ORourkeTouri2016}.
But its rank modulo two is at most \(\lceil n/2\rceil\) for
\emph{every} simple graph (Proposition~\ref{modtwo:adjacency-half-rank}).
For \(n>1\), adjacency powers can therefore never generate
\(\mathbb Z^n\). Typically, adding the degree operator removes this obstruction
together with every other prime obstruction.

\begin{theorem}[Random integral AD-saturation]
\label{integral:random}
Let \(G\sim G(n,1/2)\). Then
\[
 \Pr\bigl(\mathbb Z\langle A_G,D_G\rangle e\neq\mathbb Z^n\bigr)
       \leq n^{-3/20+o(1)}.
\]
In particular, \(\Lambda_G=\mathbb Z^n\) with probability tending to one.
\end{theorem}

On the AD-full locus, the theorem asserts that the entire finite defect group
\(C_G=\mathbb Z^n/\Lambda_G\) vanishes, not just its component at a
prescribed prime. Indeed, Smith normal form gives the equivalence
\begin{equation}\label{paper:localglobal}
 \Lambda_G=\mathbb Z^n
 \quad\Longleftrightarrow\quad
 \mathbb F_p\langle A\bmod p,D\bmod p\rangle(e\bmod p)
       =\mathbb F_p^n\quad\text{for every prime }p.
\end{equation}
The order of these quantifiers matters: all fields are handled on
the same random-graph event. The second generator is determined by
the first, so the theorem is also a statement about how a local
operation derived from an operator changes its integral reachability.

This arithmetic form of generation provides more information than
a basis over \(\mathbb Q\). It eliminates denominators in the
expressions for all vertex indicators, makes the same seed cyclic
in every characteristic, and gives a concrete obstruction theory
for reconstructing vertex coordinates from scalar observations.
The proof has a corresponding structure: construct the coordinates
with controlled denominators, then remove the finitely many
possible prime factors. This separates a global integer construction
from local invariant-subspace estimates. The exponent \(3/20\)
comes from degree thresholds and is not asserted to be optimal.

\paragraph{The construction and its consequences.}
Write \(\mathcal D_G\) for the distinct degrees, \(q=|\mathcal D_G|\),
and
\[
 \Delta_G=\prod_{\substack{d<h\\d,h\in\mathcal D_G}}(h-d).
\]
The main constructive step proves, with the same failure exponent,
\begin{equation}\label{intro:integer-inclusion}
 \Delta_G^3\mathbb Z^n\subseteq
 \operatorname{span}_{\mathbb Z}
 \{D^ie,\ D^iAD^je,\ D^iAD^jAD^\ell e:0\le i,j,\ell<q\}.
\end{equation}
Integer-surjective adjacency cuts propagate coordinates from
singleton degree classes through two rounds of recovery.
This handles every prime not dividing \(\Delta_G\) at once.
All remaining primes are at most the degree width, hence at most
\(n-1\); characteristic two, fixed odd primes, and growing primes
require different treatments.

For a word \(w\), the \emph{principal moment}
\(\mu_G(w)=e^{\mathsf T}w(A,D)e\) is a scalar observation of this
representation. On the AD-full locus these moments determine the
pointed representation up to rational orthogonal equivalence.
Integral saturation makes that equivalence a permutation.
The integer construction yields the following more economical
set of observations.

\begin{theorem}[Small-diameter identification; Corollary~\ref{cor:random-caterpillars}]
\label{paper:smalltests}
Let \(\mathcal F_n\) be the family of all ordinary caterpillars with diameter
at most seven and at most
\(\lceil12\sqrt{n\log n}+6\rceil\) vertices, including \(K_1\).
It has \(O(n^3(\log n)^3)\) isomorphism types. With probability at
least \(1-n^{-3/20+o(1)}\), \(G\sim G(n,1/2)\) is determined among
all finite simple graphs by
\((\operatorname{hom}(T,G))_{T\in\mathcal F_n}\).
\end{theorem}

The family is fixed before \(G\) is drawn, its tests are unlabelled
and unweighted, and each test has sublinear order.
This is an identification theorem: no polynomial-time decoder
from this profile is claimed. Its proof uses the integer cut
construction but does not need the local-prime estimates.
Those estimates establish the stronger integral-generation theorem.
At the other extreme, even all caterpillars of diameter at most
four fail to identify almost every graph. This qualitative lower
bound also follows from the switching argument in
Pikhurko--Verbitsky~\cite[Lemma~6.2]{PikhurkoVerbitsky2011};
Proposition~\ref{prop:caterpillar-four} supplies an entropy proof.
The gap between diameters five and seven remains open.

\paragraph{Related work.}
Our starting point is the difference between rational controllability
and integer generation. O'Rourke--Touri~\cite{ORourkeTouri2016}
establish the former for random graphs. Van Werde~\cite{VanWerde2025}
studies walk cokernels for independent-entry random matrices,
whereas here symmetry and \(D=\operatorname{diag}(Ae)\) impose
dependencies. Wood's random-sandpile work~\cite{Wood2017}
provides a symmetric random-matrix estimate used directly in our
fixed-prime argument. Smith invariants also have an algorithmic
role in exact linear algebra~\cite{KannanBachem1979,EberlyGiesbrechtVillard2000}.
Here the issue is the typical disappearance of all invariant
factors for a word-generated lattice.

The passage from algebraic obstructions to algorithmic structure
has a broader precedent in noncommutative rank and operator
scaling~\cite{GargGurvitsOliveiraWigderson2016,AllenZhuGargLiOliveiraWigderson2018}.
Those works concern rank and invariant-theoretic obstructions
over fields; they do not imply integer saturation of a cyclic
module. Our local invariant-subspace arguments and global
denominator construction address the additional arithmetic problem.

For graph reconstruction, walk-lattice and denominator methods appear
in~\cite{Wang2017,WangWangZhu2023}; unimodular overlattices also
appear in Xiang's recent preprint~\cite[Section~5]{Xiang2026}.
That preprint announces determination of almost all graphs by their
generalized spectrum. Our proof is independent, and generalized
cospectrality does not require intertwining \(D\).
The search for finite adjacency--degree descriptions of trees
appeared earlier in the setting of hyperbolic set
comparison~\cite[Sections D.2 and E.1]{LiEtAl2025}.
Lu--Li~\cite{LuLi2026} develop principal AD moments and their
caterpillar interpretation, while~\cite{Lu2026Probes} constructs
uniform probes for real cyclic modules. The distinction here is
integral generation and its random-graph saturation.
The caterpillar/degree-coloured-walk correspondence also appears
in~\cite{Cerny2025}; all tree homomorphism counts characterize
colour-refinement equivalence~\cite{DellGroheRattan2018}.
Neither correspondence supplies the fixed diameter and sublinear
test order of Theorem~\ref{paper:smalltests}.
Degree information also enters spectral reconstruction through matrix
pencils and multivariate spectra~\cite{WangLiLuXu2011,LiSun2021,WangZhao2026};
simultaneous similarity of \(A,D\) is studied
in~\cite{GodsilSun2025,FanXingZhangWang2025,ZhangFan2026}.
Bounds for generalized cospectral fibres appear
in~\cite{RazaShabbirAbbas2025,RazaAhmadShabbirAbbas2026}.
The integer word lattice and its random saturation are the additional
objects studied here.
Our coordinate propagation uses the separated-degree mechanism
of Babai--Erd\H{o}s--Selkow~\cite{BabaiErdosSelkow1980}, with integer
surjectivity replacing distinct combinatorial signatures.

\paragraph{Outline.}
The proof follows the arithmetic obstruction rather than the
list of prime estimates. Section~\ref{paper:construction} first
produces vertex coordinates with a common explicit denominator,
reducing an unbounded problem to a bounded prime range.
Section~\ref{paper:primes} then removes that denominator prime by
prime while respecting the dependence between adjacency and degrees.
Section~\ref{paper:reconstruction} explains why integer coordinates
force graph rigidity, how the short construction yields bounded-diameter
tests, and what remains between identification and decoding.
\ifdefined\PaperSplitBundle
Appendices~\ref{sec:foundations}--\ref{sec:caterpillar-reconstruction}
give the finite framework, the full integer construction and final
saturation argument, and the caterpillar proof. The separate supplement
contains the full local-prime proofs, nonsaturated fibre theory,
index-two classification, and exact computational certificates.
\else
The appendices give all complete proofs, including the local-prime
estimates, nonsaturated fibre theory, index-two classification,
and exact computational certificates.
\fi

%% file: integer_construction.tex
\section{Generating coordinates with controlled denominators}
\label{paper:construction}

\subsection{The arithmetic obstructions}

To see exactly what must be proved, let \(r=\dim_{\mathbb Q}\mathbb Q\Lambda_G\)
and choose any finite integer generating matrix for \(\Lambda_G\).
Such a matrix exists because subgroups of \(\mathbb Z^n\) are finitely
generated. Its nonzero Smith entries \(s_1\mid\cdots\mid s_r\) give
\[
 \mathbb Z^n/\Lambda_G
 \simeq \mathbb Z^{n-r}\oplus
       \bigoplus_{i=1}^r\mathbb Z/s_i\mathbb Z.
\]
The free part records missing rational directions; the torsion
part records divisibility within the visible space.
Reduction modulo \(p\) is surjective exactly when
\(r=n\) and no \(s_i\) is divisible by \(p\).
This proves \eqref{paper:localglobal}.
In particular, rational generation alone leaves every nontrivial
Smith entry unresolved.

The obstruction for adjacency alone is especially concrete.
Over \(\mathbb F_2\), put \(m_j=e^{\mathsf T}A^je\).
For \(j\ge0\), alternation gives
\[
 m_{2j+1}=(A^je)^{\mathsf T}A(A^je)=0,
 \qquad
 m_{2j}=(A^je)^{\mathsf T}(A^je)=m_j.
\]
Repeatedly dividing a positive even index by two shows \(m_j=0\)
for all \(j>0\). Consequently all positive-degree walk vectors
are pairwise orthogonal and orthogonal to \(e\).
For even \(n\), the entire walk module is totally isotropic,
of dimension at most \(n/2\).
For odd \(n\), the positive-degree part is totally isotropic in
the nondegenerate symplectic space \(e^\perp\), of dimension \(n-1\);
adding \(e\) gives the bound \(\lceil n/2\rceil\).
Thus this is an obstruction present in every graph, rather than
a rare failure of random rational controllability.

For each prime write
\[ M_p(G)=\mathbb F_p\langle A\bmod p,D\bmod p\rangle(e\bmod p). \]
For two symmetric operators, the relevant obstruction has an exact
linear-algebraic form. Over any field, the orthogonal complement of
\(M=\langle A,D\rangle e\) is invariant under \(A,D\), and \(M\)
is proper precisely when there is a nonzero common invariant
subspace \(W\subseteq e^\perp\).
Indeed, symmetry implies
\(\langle Ax,w\rangle=\langle x,Aw\rangle\), and similarly for \(D\).
Conversely every word vector is orthogonal to such a \(W\).
The ambient dot product is nondegenerate in every characteristic,
even when its restriction to \(W\) is degenerate.
The dimensions of these restricted radicals will matter in the
local counts. If all degrees are congruent to a scalar \(c\)
modulo \(p\), then \(D=cI\) and \(Ae=ce\), so
\(M_p=\mathbb F_pe\); no unconditional saturation statement can hold
for all graphs.

There is also a useful finite-length certificate. For a graph on
\(n\) vertices,
\begin{equation}\label{paper:short-integral-certificate}
 \Lambda_G=\mathbb Z^n
 \quad\Longleftrightarrow\quad
 \operatorname{span}_{\mathbb Z}\{w(A,D)e:|w|\le n-1\}=\mathbb Z^n.
\end{equation}
Over a field, the spaces \(V_j\) generated by words of length at
most \(j\) satisfy \(V_{j+1}=V_j+AV_j+DV_j\).
An equality is permanent, so a full cyclic module reaches dimension
\(n\) by step \(n-1\).
When \(\Lambda_G=\mathbb Z^n\), apply this argument over every prime
field and use the Smith criterion on the same finite integer
generating matrix. This proves the forward implication; the reverse
is immediate. The statement bounds the length of a saturation
certificate, not the stabilization length of an arbitrary word lattice.
It also does not bound the number of words by a polynomial.

\subsection{One construction for the unbounded prime range}

The first difficulty in \eqref{paper:localglobal} is the unbounded
set of primes. Proving a high-probability statement for each fixed
prime would not resolve it. Instead we construct a sublattice of
explicit finite index using words with only two occurrences of \(A\).

For a degree \(d\in\mathcal D_G\), let \(P_d\) be the diagonal projector
onto its degree class. Lagrange interpolation gives
\begin{equation}\label{paper:projectors}
 P_d=\prod_{\substack{h\in\mathcal D_G\\h\ne d}}\frac{D-hI}{d-h},
 \qquad \Delta_G P_d\in\mathbb Z[D].
\end{equation}
Each cleared projector has degree at most \(q-1\).
Consequently a singleton degree vertex \(v\) is accessible as
\(\Delta_G e_v=\Delta_GP_de\).
To reach nonsingleton classes, the construction needs integer
surjectivity, not merely full rank of a cut matrix.

\begin{lemma}[Integer-surjective cuts; Lemma~\ref{largeprime:rectangular}]
For a \(k\times m\) matrix \(B\) of independent uniform
\(\{0,1\}\) entries, \(m\ge16k\), and \(h=\lfloor m/2\rfloor\),
\[
 \Pr[B\mathbb Z^m\ne\mathbb Z^k]
 \le\left(1+\frac k2\log_2\max(2,k)\right)e^{-h/16}.
\]
\end{lemma}

Here is the reason a finite random matrix can control every prime.
Over any field, a fair Bernoulli column lies in a fixed proper
subspace with probability at most \(1/2\): condition on all but a
coordinate with a nonzero coefficient in a separating linear form.
The first \(h\) columns therefore have rational rank \(k\) except
with probability \(e^{-h/16}\). Select a nonzero \(k\times k\)
minor \(b\) from them. Hadamard gives \(|b|\le k^{k/2}\), so only
\((k/2)\log_2\max(2,k)\) primes can divide it. Every other prime is
already settled. Conditional on the first half, the second half
is independent and supplies full rank at each prime on this
finite list with the same exponential failure bound. Smith
normal form finishes the argument.
This two-stage exposure avoids a union bound over infinitely many primes.

\paragraph{Why two propagation rounds suffice.}
The random graph supplies three scales of degree classes:
many vertices with isolated extreme degrees, a large population
in small low-degree classes, and arbitrary classes of size
\(O(\sqrt n)\). The parameters in the full proof are
\[
 k=n^{0.02+o(1)},\quad r=n^{0.14+o(1)},\quad
 K=10\sqrt n+O(1),\quad s=n^{0.505+o(1)}.
\]
There are \(r\) usable extreme-degree seeds, each in a singleton
class. Low-degree classes have fewer than \(k\) vertices, while
there are at least \(2s\) vertices below a slightly lower degree
threshold. All degree classes have fewer than \(K\) vertices.
The useful inequalities are
\[
 k\log n=o(r),\qquad K\log n=o(s).
\]
They allow the exponentially small cut failure probabilities to
be union-bounded over all possible target sets.

If coordinates on a seed set \(S\) are already available and
\(C\) is a target degree class with
\[
                 A[C,S]\mathbb Z^S=\mathbb Z^C,
\]
then every \(e_u\), \(u\in C\), is an integer combination of
\(P_CA e_v\), \(v\in S\).
Using \eqref{paper:projectors} gives the following progression:
\[
 \Delta_G e_v\ \text{at singleton seeds}
 \ \longrightarrow\
 \Delta_G^2e_u\ \text{in small low-degree classes}
 \ \longrightarrow\
 \Delta_G^3e_u\ \text{at every vertex}.
\]
The first round uses \(r\gg k\) columns. The second uses
\(s\gg K\) columns among the recovered low-degree vertices.
Each round adds one \(A\) and one cleared degree projector.
The resulting vectors therefore lie in the span of
\(D^ie,D^iAD^je,D^iAD^jAD^\ell e\), with \(0\le i,j,\ell<q\).

\paragraph{Preserving independence.}
Conditioning on the final degree classes and then treating their
cut edges as independent would be invalid. The proof reverses
that order. Fix a target row set \(C\), expose the graph outside
\(C\), and choose candidate columns from this exposed graph.
The edges between \(C\) and these columns are still independent
Bernoulli variables. Prove integer surjectivity simultaneously
for every eligible \(C\), then substitute the actual degree
classes of the full graph.

For the first round, the largest \(r+1\) degrees are separated by more
than \(g=n^{0.03+o(1)}\). Since \(k=o(g)\), deleting a target set \(C\) disjoint from the seed set \(S\),
with \(|C|\le k\), leaves the seed choice stable. A degree class
meeting \(S\) is already a recovered singleton.
For the second round, the selected columns have degree at most
\(t-K\) in the exposed graph, where
\[
 t=\left\lfloor (n-1)/2-\sqrt{0.245\,n\log n}\right\rfloor.
\]
There are at least \(2s-K\ge s\) eligible vertices outside \(C\).
Adding back at most \(K\) vertices keeps their full degrees at most
\(t\), so they belong to classes already recovered in the first
round. The uniform cut statements can now be applied without
conditioning on the cut entries themselves.

The degree separation and count estimates, together with these
uniform cut events, fail with probability at most
\(n^{-3/20+o(1)}\). This proves the integer inclusion
\eqref{intro:integer-inclusion}; the complete argument is in
Appendix~\ref{largeprime:section}. It has two consequences to
carry forward. First, \(\Lambda_G\) has full rational rank and
the exponent of \(\mathbb Z^n/\Lambda_G\) divides \(\Delta_G^3\).
Second, every possible bad prime divides a difference of graph
degrees and hence is at most \(n-1\). What remains is to remove
these prime obstructions without assuming independence between
\(A\) and its degree operator.

%% file: prime_mechanisms.tex
\section{Local mechanisms for integral generation}
\label{paper:local-mechanisms}
\label{paper:primes}

The integer inclusion confines every possible prime divisor of the
defect to the degree differences. We must still exclude these primes
without treating the degree diagonal as independent of the adjacency
matrix. The common obstruction is a nonzero subspace
\(W\subseteq e^\perp\) invariant under both \(A\) and \(D\): symmetry
identifies such a space with a failure of cyclicity. The way we count
these obstructions changes with the characteristic. Full proofs are
in \ifdefined\PaperSplitBundle Supplement Sections\else Appendices\fi~\ref{modtwo:section}, \ref{oddprime:section}, and
\ref{sec:random-growing-primes}.

\subsection{Parity: count inside the actual degree fibre}

\begin{theorem*}[Prime two; Theorem~\ref{modtwo:random}]
Let $G\sim G(n,1/2)$. For every $n\geq 48$,
\[
 \Pr\bigl[M_2(G)\neq\mathbb F_2^n\bigr]
 \leq
 4\exp(-n/18)+1024(n+1)^4\,2^{-n/48}.
\]
Consequently $M_2(G)=\mathbb F_2^n$ with probability
$1-\exp(-\Omega(n))$.
\end{theorem*}

Over \(\mathbb F_2\), adjacency matrices are uniform alternating
matrices. The map \(A\mapsto Ae\) is onto \(e^\perp\), since a single
edge maps to the sum of its endpoint vectors. Thus \(d=Ae\) is uniform
on \(e^\perp\), and conditioning on \(d\) leaves a uniform affine
space, not an independent-edge model. Except with probability
\(4\exp(-n/18)\), both parity classes \(V_0,V_1\) have sizes
\(n_0,n_1\geq n/3\).

Fix such a balanced fibre. Let \(e_i\) be the all-ones vector on
\(V_i\). Since \(D\) projects onto \(V_1\), an
invisible invariant subspace decomposes as
\(W=W_0\oplus W_1\), with \(W_i\subseteq e_i^\perp\).
Write \(t_i=\dim W_i\), let \(r_i\) be the radical dimension of its restricted dot product, and put
\(\varepsilon_i=\mathbf1_{\{e_i\in W_i\}}\),
\(t=t_0+t_1\), and \(r=r_0+r_1\).
The invariance constraints are the bivectors
\(W\wedge W^\perp\), whereas the prescribed degree vector fixes
\(e\wedge\mathbb F_2^n\). Their dimensions and intersection give
\[
 c=t(n-t)-\frac{r(r+1)}2,\qquad
 c_{\rm cond}=
 \begin{cases}
 c-t,&e\notin W,\\
 c-n+r+1,&e\in W.
 \end{cases}
\]
For these graded invisible spaces the compatibility conditions are
automatic. The conditional probability of invariance is exactly
\(2^{-c_{\rm cond}}\).

The count must retain the same radical information. In a block of
size \(n_i\), the even-weight hyperplane is symplectic when \(n_i\)
is odd and has radical \(\langle e_i\rangle\) when \(n_i\) is even.
Counting through its symplectic quotient gives at most
\(\gamma^4 2^{E_0+E_1}\) candidates in a fixed parameter stratum, where
\[
 \gamma=\prod_{j\geq1}(1-2^{-j})^{-1}<4,\qquad
 E_i=t_i(n_i-1-t_i)-\frac{r_i(r_i-1)}2
       -\varepsilon_i(n_i-t_i-r_i).
\]
The last term is needed near full dimension. Set
\(X=t_0(n_1-t_1)+t_1(n_0-t_0)\) and
\(m=\min(t,n-t)\). The exact formulas imply
\[
 c_{\rm cond}-(E_0+E_1)\geq X/2-r
 \geq (n/24-1)m\geq n/48.
\]
Here balance gives \(X\geq nm/12\), and \(r\leq m\).
Summing over at most \(4(n+1)^4\) parameter strata proves the
conditional bound uniformly in \(d\); averaging over the fibres
proves the theorem. This argument excludes invisible subspaces of
all dimensions, including those larger than \(n/2\).

\subsection{Fixed odd primes: keep the seed in the code}

\begin{theorem*}[Fixed odd primes; Theorem~\ref{odd:fixed-irreducible}]
For every fixed odd prime $p$, there are constants $c_p,K_p>0$ such that
\[
 \Pr\bigl[\text{$A,D$ have a nonzero proper common invariant subspace}\bigr]
 \leq \frac{K_p}{n}+K_p\exp(-c_p n).
\]
Consequently
$\mathbb F_p\langle A,D\rangle e=\mathbb F_p^n$ with probability tending
to one.
\end{theorem*}

The conclusion is irreducibility over \(\mathbb F_p\), which is
stronger than the required cyclicity. Symmetry reduces the possible
dimensions to \(t\leq n/2\). Degree-residue classes have sizes
\(n/p\pm\epsilon n\) except with probability \(O_{p,\epsilon}(n^{-1})\).
We enumerate deterministic balanced assignments \(g\), and estimate
each event under the original edge distribution.

For a fixed \(t\)-space with radical dimension \(r\), invariance
has codimension
\(c_{\rm sym}=t(n-t)-r(r-1)/2\) in the symmetric matrices.
Removing the diagonal costs at most \(n\) constraints, so its
Bernoulli probability is at most \(2^{-c_{\rm sym}+n}\).
Counting degree-graded subspaces together with their block radicals
then gives, for sufficiently small fixed \(\epsilon\), an exponent
\[
 -\theta nt/8+(\alpha+1)n,
 \qquad \alpha=\log_2p,\quad \theta=1-\alpha/p>0.
\]
Thus dimensions \(t\geq T\), for a sufficiently large constant
\(T=T(p)\), are excluded exponentially. The remaining bounded
dimensions require more than this coarse constraint bound.

Let the rows of \(U\) span a small candidate \(W\). If \(e\notin W\),
adjoin \(e^{\mathsf T}\) to form \(F=[U;e^{\mathsf T}]\); if
\(e\in W\), choose a basis matrix \(F\) whose last row is
\(e^{\mathsf T}\). In either case every column has final coordinate
one. Degree grading puts its other coordinates in prescribed
eigenspaces. The equations for adjacency invariance, together with
\(e^{\mathsf T}A=g^{\mathsf T}\), prescribe a target for \(FA\).
Keeping this seed row is what makes the target estimate compensate
for enumerating the \(p^n\) degree assignments.

A full-row-rank matrix is a code of distance \(\delta n\) if deleting
fewer than \(\delta n\) columns preserves its rank. For a code with
\(k\) rows, Wood's arbitrary-target bound
\cite[Lemma~4.1]{Wood2017}, applied after adding an independent
Bernoulli diagonal, gives
\[
 \Pr[FA=Y]\leq 2^n Kp^{-kn}.
\]
The factor \(2^n\) is the cost of conditioning that diagonal to zero.
Degree grading allows at most \(p^{(k-1)(n/p+\epsilon n)}\) choices
of \(F\) for each assignment and fixed auxiliary matrices. Including
the assignments gives
\[
 2^n p^{-(k-1)(1-1/p-\epsilon)n}.
\]
This decreases exponentially: \(k\geq2\), and
\((1-1/p)\log p>\log2\) for every odd prime.

For a non-code, choose a maximal depth witness: deleting a small set
\(\tau\), of size \(s\), leaves columns spanning a proper space
\(H\) of dimension \(h\), and the remaining matrix \(F_\eta\) is a
code, where \(\eta=[n]\setminus\tau\). Maximality supplies the code
property. The equations projected modulo \(H\) use only cut edges
between \(\tau\) and \(\eta\). A nonzero quotient coefficient forces
their probability to be at most \(2^{-|\eta|}\). After those edges
are exposed, the internal matrix on \(\eta\) is independent, so
the code estimate gives
\[
 \Pr[FA=Y]\leq
 2^{-|\eta|}\,2^{|\eta|}Kp^{-h|\eta|}
 =Kp^{-h(n-s)}.
\]
The quotient equations pay for the zero-diagonal cost. For \(h\geq2\),
the degree-graded affine slices have total direction dimension at
most \(h-1\); counting them makes this estimate summable over all
small exceptional sets. If \(h=1\), every surviving column is the
same. Every residue class still occurs, forcing \(W\) to be
supported on \(\tau\) (and making the encoding with \(e\in W\)
impossible). Independent outside rows then give
\(2^{-t(n-s)}\), which dominates the number of small-support spaces.
This completes the bounded-dimensional part without assuming that
the degree partition is independent of the edges.

\subsection{Growing primes: split by support and dimension}

\begin{theorem*}[Growing primes; Theorem~\ref{growing:uniform}]
There is an absolute integer \(P_0\) such that, for
\(G\sim G(n,1/2)\), with probability \(1-\exp(-\Omega(\sqrt n))\) the pair
\((A_G,D_G)\) has no nonzero proper common invariant subspace over
\(\mathbb F_p\), simultaneously for every prime \(P_0\le p\le n\).
\end{theorem*}

The fixed-prime constants are not used in this growing interval.
Uniformly for integer moduli \(128\leq p\leq n\), every degree-residue
class has size less than
\[
 b_p=\left\lceil32\left(\frac np+\sqrt n\right)\right\rceil
\]
except on one event of probability \(\exp(-\Omega(\sqrt n))\).
For a deterministic proposed class, expose its internal edges.
The external degrees are independent binomials; unimodality bounds
each residue mass by \(1/p+4/\sqrt n\). A union bound over the
proposed classes and all moduli gives the simultaneous statement.

Put \(\ell_p=\log_2p\) and \(T_p=\lceil16(\ell_p+1)\rceil\).
Choose the absolute \(P_0\) first and then take \(n\) sufficiently
large that, uniformly for \(P_0\leq p\leq n\),
\[
 T_pb_p\leq\frac{n}{100\ell_p},\qquad
 b_p\ell_p\leq n/32,\qquad \log_2(4n)\leq n/32.
\]
This order is possible because
\(\sup_{p\geq P_0}(\log p)^2/p\to0\) and
\((\log n)^2/\sqrt n\to0\).

A degree-graded space of dimension \(t<T_p\) is supported on at
most \(tb_p\) coordinates. Any nonzero vector in it is annihilated
by the outside cut. Taking a union over all supports and vectors,
without enumerating degree assignments, bounds this probability by
\[
 \sum_{1\leq j\leq n/(100\ell_p)}
       \binom nj p^j2^{-(n-j)}
 \leq(n+1)2^{-0.89n}.
\]
For larger dimensions, symmetry lets us replace an invariant space
by its orthogonal complement, so it suffices to consider
\(T_p\leq t\leq n/2\). Enumerate at most \(p^n\) assignments and
\((4n)^tp^{b_pt}\) graded subspaces for each assignment.
The joint event that an assignment is the actual one and that
\(AW\subseteq W\) is contained in the latter event alone.
The unconditional invariance bound therefore applies, giving binary
exponent at most
\[
 n\ell_p+t\log_2(4n)+b_pt\ell_p-nt/4+n
 \leq -2n(\ell_p+1).
\]
Summing over dimensions and at most \(n\) primes preserves the
claimed failure bound. This is uniform control of an entire interval,
not a union of estimates whose constants depend on a growing prime.

The absolute cutoff leaves only finitely many fixed odd primes.
Combining their bounds with the parity theorem, the growing interval,
and the integer inclusion gives
\[
 \Pr[\Lambda_G\ne\mathbb Z^n]
 \leq n^{-3/20+o(1)}+\exp(-\Omega(n))
        +O(n^{-1})+\exp(-\Omega(\sqrt n))
 =n^{-3/20+o(1)}.
\]
The local criterion applies on this single intersection event and
eliminates every Smith factor simultaneously.

%% file: reconstruction.tex
\section{From integer generation to graph reconstruction}
\label{paper:reconstruction}

The arithmetic theorem also explains when scalar observations
recover a graph's vertex coordinates. For noncommuting words
\(w\) in two letters, put
\[
 \mu_G(w)=e^{\mathsf T}w(A,D)e,\qquad
 w^*(A,D)=w(A,D)^{\mathsf T}.
\]
The moments give inner products of word vectors:
\(\langle u(A,D)e,v(A,D)e\rangle=\mu_G(u^*v)\).
Thus they determine the Euclidean representation generated by \(e\).
The remaining issue is its integer coordinate lattice.

\paragraph{The deterministic bridge.}
If \(G\) is AD-full of order \(n\), choose a word basis
\[ B_G=[w_1(G)e,\ldots,w_n(G)e], \]
with \(w_1\) empty and
\(|w_i|\le n-1\). Matching moments through length \(2n-1\) give
\[
 B_G^{\mathsf T}B_G=B_H^{\mathsf T}B_H,\qquad
 B_G^{\mathsf T}X_GB_G=B_H^{\mathsf T}X_HB_H
 \quad(X=A,D).
\]
The empty-word moment fixes the order of \(H\). The Gram identity
makes \(B_H\) invertible, and the unique
\(Q=B_HB_G^{-1}\) satisfies
\[
 Q^{\mathsf T}Q=I,\quad Qe=e,\quad
 QA_G=A_HQ,\quad QD_G=D_HQ,\quad Q\Lambda_G=\Lambda_H.
\]
Conversely these identities give every moment equality.
This proves the finite pointed reconstruction theorem
(Theorem~\ref{thm:finiteGNS}).

If \(\Lambda_G=\mathbb Z^n\), the last identity makes every column
of \(Q\) integral. An integral unit vector is a signed coordinate
vector; orthogonality makes the supports distinct, and \(Qe=e\)
makes the signs positive. Hence \(Q\) is a permutation matrix,
and \(G\cong H\). Only \(G\) needs to satisfy the saturation hypothesis.

\paragraph{Why short words give short test trees.}
Let \(C(a_0,\ldots,a_k)\) be a path with \(k\) edges and \(a_i\)
additional leaves at its \(i\)th vertex. Summing the images of
its leaves gives the established dictionary~\cite{LuLi2026,Cerny2025}
\[
 \operatorname{hom}(C(a_0,\ldots,a_k),G)
       =e^{\mathsf T}D^{a_0}AD^{a_1}\cdots AD^{a_k}e.
\]
A word with \(k\) occurrences of \(A\) and length \(L\) therefore
uses a tree of diameter at most \(k+2\) and order \(L+1\).

On the integer-construction event, take \(B_G\) from the words
in \eqref{intro:integer-inclusion}. Each basis word has at most
two \(A\)'s and length at most \(3q-1\). The Gram and inserted
operator matrices above use at most five \(A\)'s and length
at most \(6q-1\). All required counts thus come from caterpillars
of diameter at most seven and order at most \(6q\).
Degree concentration gives
\(q\le2\sqrt{n\log n}+1\), except with probability \(2/n\).
The number of caterpillar types of order at most \(m\) and diameter
at most seven is at most
\(\sum_{j=1}^{6}\binom mj=O(m^6)\).
These observations give the quantitative bounds of
Theorem~\ref{paper:smalltests}.

There is one further point in that argument: to turn this short
Gram reconstruction into a permutation, one need not invoke the
local-prime estimates. An orthogonal pointed AD intertwiner already
maps singleton degree vertices to corresponding positive unit
vectors. Suppose it has done so for a set \(S\), and
\(A_G[C,S]\mathbb Z^S=\mathbb Z^C\). Each \(e_u\) in \(C\) is an
integer combination of \(P_CA_Ge_v\), \(v\in S\). Its image is
therefore an integer combination of the integer vectors
\(P_{C'}A_He_{v'}\). Orthogonality and pointedness again make
it a positive coordinate unit vector. The two cut rounds
propagate this conclusion to every vertex
(Remark~\ref{rem:implications-shorter-rigidity}).

This uses the actual integer-surjective cuts, not just the inclusion
\(\Delta_G^3\mathbb Z^n\subseteq\Lambda_G\).
It explains the logical division of the results: the global
construction already gives typical moment rigidity, whereas
the local analysis proves the stronger assertion that the original
word lattice itself has no arithmetic defect.

\paragraph{What a nonzero defect permits.}
For an AD-full graph, let
\(\Lambda^\#=\{x\in\mathbb Q^n:\langle x,\Lambda_G\rangle
\subseteq\mathbb Z\}\). Theorem~\ref{thm:arith-fibre} identifies
the graph isomorphism classes with moments \(\mu_G\) with the
\(A,D\)-stable integral unimodular lattices
\[
                    \Lambda_G\subseteq L\subseteq\Lambda^\#
\]
that admit an orthonormal basis \(v_1,\ldots,v_n\) satisfying
\(e=\sum_i v_i\), in which \(D\) is diagonal and \(A\) is a
zero-diagonal symmetric \(\{0,1\}\) matrix.
These are exact graph-realization conditions.
The identity \(Ae=De\) then ensures that \(D\) records the degrees.
Each admissible lattice gives exactly one graph isomorphism class.

This description is finite but is not itself a polynomial-time
classification algorithm: finding a suitable unit frame remains
part of admissibility. At index two it becomes explicit.
Let \(v\) be the unique nonzero binary vector orthogonal to all
word vectors and \(S=\operatorname{supp}(v)\).
The fibre has at most two classes. It has two exactly when
\(|S|=4\), \(G[S]\) is regular, every outside vertex has
zero, two, or four neighbours in \(S\), and the four vertices
of \(S\) have the same total degree. The other graph is obtained
by the Godsil--McKay switch on \(S\)
(Theorem~\ref{thm:arith-index-two-classification}).
We give both a rigid nonsaturated example and an
example attaining the bound of two.

\paragraph{Computation and the remaining decoding problem.}
When \(G\) is given, one can compute \(\Lambda_G\) and its Smith
invariants in deterministic polynomial time by repeatedly replacing
a Hermite basis \(B\) by one for \([B\mid AB\mid DB]\).
After rational rank stabilizes, each strict enlargement reduces
the finite saturation index by a factor of at least two.
Hadamard bounds the initial index and hence the number of rounds;
\ifdefined\PaperSplitBundle Supplement Section\else Appendix\fi~\ref{app:integral-computation} gives the bit bounds.
Saturation also has a shorter certificate: words of length at most
\(n-1\) generate \(\mathbb Z^n\) whenever all words do
(Proposition~\ref{integral:short-certificate}).

More precisely, with \(r=\dim_{\mathbb Q}\mathbb Q\Lambda_G\),
\(\Delta_*=\max(1,\max_v\deg(v))\), and
\[
 I_0=[\,\mathbb Q\Lambda_G\cap\mathbb Z^n:\Lambda_{r-1}\,],
 \qquad
 \Lambda_j=\operatorname{span}_{\mathbb Z}
                \{w(A,D)e:|w|\le j\},
\]
Theorem~\ref{thm:finiteclosure} gives
\[
 \Lambda_G=\Lambda_{r-1+\lfloor\log_2 I_0\rfloor},
 \qquad I_0\le r^{r/2}\Delta_*^{r(r-1)/2}.
\]
Choose rationally independent words in nondecreasing length.
The \(i\)th has length at most \(i-1\), since an earlier
stabilization would be permanent. Their entries have absolute
value at most \(\Delta_*^{i-1}\), giving the displayed index bound
by Hadamard. Thereafter every strict lattice enlargement divides
the remaining index by an integer at least two.
Each Hermite step uses at most \(3n\) columns, and \ifdefined\PaperSplitBundle the supplement\else the appendix\fi{} bounds their bit lengths by a polynomial in \(n\).
Thus the algorithm does not enumerate the exponentially many
words in \eqref{paper:short-integral-certificate}.

Starting only with the count profile is a different input problem.
Gram matrices and operator insertions recover an abstract rational
representation, but an arbitrary lattice basis need not be its
orthonormal vertex basis. The proof of uniqueness does not by
itself find that basis in polynomial time. Exploiting the special
random cut structure to do so is a further algorithmic question.
Likewise, diameter four is insufficient, while the present
construction reaches seven; in this construction, improving the upper bound requires
a spanning frame with fewer adjacency occurrences, not merely a smaller interpolation
denominator.

%% file: article_foundations.tex
\section{Principal modules and integral saturation}\label{sec:foundations}

All graphs in this paper are finite, simple, and undirected, with at least
one vertex.  Write \(A=A_G\), \(D=D_G=\operatorname{diag}(A_Ge)\), and
\(e=(1,\ldots,1)^{\mathsf T}\); a subscript on \(e\) specifies the graph
when needed.  A word means a word in the two
noncommuting letters \(a,d\); the empty word evaluates to the identity.
Its length is denoted by \(|w|\), and \(w^*\) is the reversed word.
Since \(A,D\) are symmetric, \(w(A,D)^{\mathsf T}=w^*(A,D)\).
Define
\[
 M_{\mathbb Q}(G)=\mathbb Q\langle A,D\rangle e,\qquad
 \Lambda_G=\mathbb Z\langle A,D\rangle e
       =\sum_w\mathbb Z\,w(A,D)e,
 \qquad r_G=\dim_{\mathbb Q}M_{\mathbb Q}(G).
\]
Here the integral algebra is unital.  Equivalently, \(\Lambda_G\) is
the smallest subgroup of \(\mathbb Z^n\) containing \(e\) and stable
under \(A,D\).  Its saturation and relative defect group are
\[
 \operatorname{Sat}(\Lambda_G)=M_{\mathbb Q}(G)\cap\mathbb Z^n,
 \qquad T_G=\operatorname{Sat}(\Lambda_G)/\Lambda_G.
\]
The group \(T_G\) is finite.  In particular,
\[
 \mathbb Z^n/\Lambda_G\simeq
       \mathbb Z^{\,n-r_G}\oplus T_G
\]
as abelian groups, without a preferred splitting.
Both operators preserve the saturation, since they preserve its
rational span and \(\mathbb Z^n\); hence they induce endomorphisms
of \(T_G\).  For a regular graph \(\Lambda_G=\mathbb Ze\) and
\(T_G=0\), so a vanishing relative defect need not mean full rank.

\begin{definition}\label{def:saturation}
The graph is \emph{AD-full} if \(r_G=n\), equivalently if
\(\operatorname{Sat}(\Lambda_G)=\mathbb Z^n\).
The lattice is \emph{primitive} (or saturated in its rational span)
if \(\Lambda_G=\operatorname{Sat}(\Lambda_G)\).
The graph is \emph{integrally AD-saturated} if
\(\Lambda_G=\mathbb Z^n\).
For a prime \(p\), put
\[
 M_p(G)=\mathbb F_p\langle\overline A,\overline D\rangle\overline e.
\]
The graph is \emph{\(p\)-AD-saturated} if \(M_p(G)=\mathbb F_p^n\).
Bars denote reduction modulo \(p\); when working over \(\mathbb F_p\)
we usually suppress them.
\end{definition}

AD-fullness and primitivity describe the two obstructions separately:
missing rational directions and integral divisibility.
For an AD-full graph write
\[
 C_G=\mathbb Z^n/\Lambda_G,\qquad
 I_G=|C_G|,\qquad c_G=\exp(C_G),
\]
with exponent \(1\) for the trivial group.

\begin{proposition}[Local criterion]\label{prop:localcriterion}
For every graph and prime \(p\), the following are equivalent:
\begin{enumerate}
\item \(G\) is \(p\)-AD-saturated;
\item \(\Lambda_G+p\mathbb Z^n=\mathbb Z^n\);
\item \(\Lambda_G\otimes_{\mathbb Z}\mathbb Z_p=\mathbb Z_p^n\);
\item \(r_G=n\) and \(p\nmid I_G\).
\end{enumerate}
Consequently integral AD-saturation is equivalent to
\(p\)-AD-saturation for every prime.
\end{proposition}
\begin{proof}
The reduction of \(\Lambda_G\) is the cyclic module over
\(\mathbb F_p\).  Take a Smith normal form of a finite
integer generating matrix for \(\Lambda_G\), with nonzero diagonal
entries \(s_1\mid\cdots\mid s_{r_G}\).  Its reduction is surjective
precisely when \(r_G=n\) and all \(s_i\) are units modulo \(p\).
The same condition makes its columns span \(\mathbb Z_p^n\).
Finally \(I_G=\prod_i s_i\) in full rank.  A finite generating matrix
exists because every subgroup of \(\mathbb Z^n\) is finitely generated.
\end{proof}

\begin{corollary}[A degree-congruence obstruction]
\label{cor:arith-degree-congruence}
If every degree is congruent to \(c\) modulo \(p\), then
\(\mathbb F_p\langle A,D\rangle e=\mathbb F_pe\).
If the graph is AD-full, then \(\dim_{\mathbb F_p}C_G/pC_G=n-1\).
\end{corollary}
\begin{proof}
Modulo \(p\), \(D=cI\) and \(Ae=ce\), so every word preserves
the line through \(e\).  The quotient statement follows by reducing
a Smith normal form.
\end{proof}

\begin{proposition}[Invariant obstruction]\label{prop:obstruction}
Over any field \(k\), the module \(k\langle A,D\rangle e\) is proper
if and only if there is a nonzero subspace
\[
             W\subseteq e^\perp,\qquad AW\subseteq W,\quad DW\subseteq W.
\]
\end{proposition}
\begin{proof}
The orthogonal complement of the cyclic module is invariant because
both operators are symmetric, and is nonzero if that module is
proper.  Conversely, invariance and symmetry imply that every vector
in \(W\) is orthogonal to every word vector \(w(A,D)e\).
The ambient dot product is nondegenerate over every field, including
characteristic two; its restriction to a subspace need not be.
\end{proof}

\begin{proposition}[Complement invariance]\label{prop:complement}
With the same vertex coordinates,
\(\Lambda_{\overline G}=\Lambda_G\).
\end{proposition}
\begin{proof}
Every subgroup of \(\mathbb Z^n\) containing \(e\) is stable under
\(J=ee^{\mathsf T}\), since \(Jx=(\sum_i x_i)e\).
Use \(A_{\overline G}=J-I-A_G\) and
\(D_{\overline G}=(n-1)I-D_G\), first in one direction and then in
the other.  No connectivity assumption is needed.
\end{proof}

\subsection{Position relative to equitable partitions}
Let \(P\) be the matrix whose columns are the indicators of the cells
of an equitable partition.  Its column space contains \(e\) and is
\(A\)-invariant.  Degrees are constant on each cell, so this space is
also \(D\)-invariant.  Therefore \(M_{\mathbb Q}(G)\) is contained in
the cell-constant space for the stable partition of color refinement.
In particular, AD-fullness forces that partition to be discrete.
It also forces asymmetry: an automorphism fixes \(e\), commutes with
\(A,D\), and hence fixes every vector in \(M_{\mathbb Q}(G)\).

Consequently color refinement identifies every AD-full graph among all
graphs: on a disjoint union, equality of the stable color data gives a
unique matching of singleton cells and their adjacencies.  The arithmetic
questions studied here concern the finer distinction between rational
span and integer generation, and reconstruction from the principal
moments defined below.  For a regular graph,
\(M_{\mathbb Q}=\mathbb Qe\), so the full-rank theory applies only when
\(n=1\).

\subsection{Finite word bounds and exact computation}\label{sec:finite}
Put
\[
 \Lambda_t=\sum_{|w|\leq t}\mathbb Z\,w(A,D)e,\qquad
 M_t=\mathbb Q\Lambda_t,\qquad
 \Delta_*=\max(1,\max_v\deg(v)).
\]

\begin{theorem}[Finite integral closure]\label{thm:finiteclosure}
Let \(r=r_G\) and
\(I_0=[M_{\mathbb Q}(G)\cap\mathbb Z^n:\Lambda_{r-1}]\).
Then \(I_0\) is finite and
\[
 \Lambda_G=\Lambda_{r-1+\lfloor\log_2 I_0\rfloor},
 \qquad
 I_0\leq r^{\,r/2}\Delta_*^{\,r(r-1)/2}.
\]
In particular, all word vectors of length at most
\begin{equation}\label{eq:wordlength}
 L=r-1+
 \left\lceil \frac r2\log_2 r+
                  \frac{r(r-1)}2\log_2\Delta_*\right\rceil
\end{equation}
generate \(\Lambda_G\).  A lattice basis and its Smith invariants can
be computed in deterministic polynomial time in the graph order.
\end{theorem}
\begin{proof}
The rational filtration reaches its final dimension by step \(r-1\).
Thereafter, every strict lattice enlargement reduces its finite
saturation index by a factor of at least two.  Hadamard's inequality
bounds the initial index, and Hermite normal form implements the
successive closures.  The details, including bit bounds, are given in
Appendix~\ref{app:integral-computation}.
\end{proof}

\input{article_short_certificate}

\subsection{Moments and their finite certificate}
The principal AD-moment function is
\[
                  \mu_G(w)=e^{\mathsf T}w(A_G,D_G)e .
\]
Write \(G\equiv_\mu H\) when these values agree for every word,
including the empty word.  The latter value records the graph order.

\begin{theorem}[Finite pointed reconstruction]\label{thm:finiteGNS}
Suppose \(G\) is AD-full of order \(n\).  The following are equivalent:
\begin{enumerate}
\item \(G\equiv_\mu H\);
\item the principal moments agree for all words of length at most
\(2n-1\);
\item there is a rational orthogonal matrix \(Q\) such that
\[
 Qe=e,\qquad QA_G=A_HQ,\qquad QD_G=D_HQ .
\]
\end{enumerate}
For fixed labelled graphs this matrix is unique, and \(H\) is AD-full.
\end{theorem}
\begin{proof}
The empty-word moment gives \(|V(H)|=n\).
Choose a word basis \(B_G=[w_1(G)e,\ldots,w_n(G)e]\), containing
the empty word, with \(|w_i|\leq n-1\), and form the corresponding
matrix \(B_H\).  Equality of the moments of length at most \(2n-2\)
gives
\[
 B_G^{\mathsf T}B_G=B_H^{\mathsf T}B_H,
\]
because the \((i,j)\) entry is \(\mu_G(w_i^*w_j)\).
Hence \(B_H\) is nonsingular and
\(Q=B_HB_G^{-1}\) is rational orthogonal and fixes \(e\).
The additional insertion identities
\[
 B_G^{\mathsf T}X_GB_G=B_H^{\mathsf T}X_HB_H,
                         \qquad X=A,D,
\]
use words of length at most \(2n-1\) and give both intertwinings.
Conversely the intertwinings transport every word vector and every
moment.  They also force \(QB_G=B_H\), proving uniqueness.
\end{proof}

Thus a recorded word basis and the three \(n\) by \(n\) matrices
\(B_G^{\mathsf T}B_G\), \(B_G^{\mathsf T}A_GB_G\), and
\(B_G^{\mathsf T}D_GB_G\) suffice to certify pointed orthogonal
equivalence.  Whether that equivalence is a vertex permutation is the
arithmetic question addressed in Section~\ref{sec:arith-intertwiners}.

%% file: article_short_certificate.tex
\begin{proposition}[A bounded-length certificate]
\label{integral:short-certificate}
For every graph on \(n\geq1\) vertices, integral AD-saturation is
equivalent to
\[
 \operatorname{span}_{\mathbb Z}
 \{w(A_G,D_G)e: |w|\leq n-1\}=\mathbb Z^n.
\]
\end{proposition}

\begin{proof}
Only the forward implication needs proof. Over any field, let
\(V_j\) be the span of word vectors of length at most \(j\).
Then \(\dim V_0=1\) and
\[
 V_{j+1}=V_j+AV_j+DV_j.
\]
Once \(V_{j+1}=V_j\), the space is invariant under both operators,
and the sequence stabilizes. Thus a full word module reaches
dimension \(n\) by step \(n-1\).

If \(\Lambda_G=\mathbb Z^n\), this argument applies modulo every
prime. The integer lattice generated by the displayed finite word
set therefore has full image modulo every prime. Its Smith
invariants are all one, so it equals \(\mathbb Z^n\).
\end{proof}

The proposition bounds the word length needed to certify full
integral generation. It does not assert that an arbitrary integer
word lattice stabilizes by length \(n-1\).

%% file: article_random_large_primes.tex
\section{Simultaneous saturation away from degree differences}
\label{largeprime:section}

We first confine the possible prime divisors of the integer defect
to a finite set determined by the degrees. Let $\mathcal D_G$ be the
set of distinct degrees of $G$, put $q=|\mathcal D_G|$, and define
\[
 \Delta_G=\prod_{\substack{d,h\in\mathcal D_G\\d<h}}(h-d),
 \qquad
 w_G=\max_{v\in[n]}\deg_G(v)-\min_{v\in[n]}\deg_G(v).
\]
An empty product is $1$. In particular, $\Delta_G$ is always a positive
integer.
Let $\Lambda_G^{[2]}$ be the integer span of
\[
 D^i e,\qquad D^iAD^j e,\qquad D^iAD^jAD^\ell e,
 \qquad 0\leq i,j,\ell\leq q-1.
\]
Every generating word has at most two occurrences of $A$ and length
at most $3q-1$.

\begin{theorem}[Integer generation with controlled denominators]
\label{largeprime:main}
If $G\sim G(n,1/2)$, then
\[
 \Pr\bigl[\Delta_G^3\mathbb Z^n
                 \not\subseteq\Lambda_G^{[2]}\bigr]
 \leq n^{-3/20+o(1)}.
\]
In particular, $\Delta_G^3\mathbb Z^n\subseteq\Lambda_G$ with
probability $1-o(1)$. On this event, $\mathbb Z^n/\Lambda_G$ is finite
and its exponent divides
$\Delta_G^3$. Consequently, on one event of probability $1-o(1)$,
\[
 \mathbb F_p\langle A\bmod p,D\bmod p\rangle(e\bmod p)
   =\mathbb F_p^n
 \quad\text{for every prime }p\nmid\Delta_G.
\]
In particular, this holds simultaneously for every prime $p>w_G$.
\end{theorem}

The proof constructs the inclusion over the integers, so the conclusion
for all $p\nmid\Delta_G$ follows on a single graph event. Degree
interpolation first produces $\Delta_G e_v$ at vertices whose degrees
are unique. Two rounds of integer linear combinations then produce
$\Delta_G^2 e_v$ throughout a large collection of small degree classes,
and finally $\Delta_G^3 e_v$ at every vertex.

The probabilistic input is integer surjectivity of rectangular adjacency
blocks. To preserve independence, we choose the columns after exposing
the graph outside a fixed row set, prove the required property
simultaneously for every eligible row set, and only then apply it to
the degree classes of the full graph. The remaining sections will
remove the possible prime divisors of $\Delta_G$.

\subsection{An integer-surjectivity estimate for rectangular matrices}

\begin{lemma}\label{largeprime:rectangular}
Let $B$ be a $k$-by-$m$ matrix with independent uniform $\{0,1\}$
entries, where $k\geq1$ and $m\geq16k$, and put $h=\lfloor m/2\rfloor$.
Then
\[
 \Pr\bigl[B\mathbb Z^m\neq\mathbb Z^k\bigr]
 \leq
 \left(1+\frac{k}{2}\log_2\max(2,k)\right)\exp(-h/16).
\]
\end{lemma}

\begin{proof}
Let $F$ be either $\mathbb Q$ or a prime field, let $U$ be a proper
linear subspace of $F^k$, and let $X$ be a uniform $\{0,1\}^k$ column.
Choose a nonzero linear functional vanishing on $U$ and a coordinate
whose coefficient in this functional is nonzero. Conditional on all
other coordinates, at most one of the two possible values of the chosen
coordinate makes the functional vanish. Therefore
\[
 \Pr[X\in U]\leq\frac12.
\]
In a sequence of independent columns, every rank increment before full
rank is reached thus has conditional probability at least $1/2$.
Coupling with independent fair bits, and using dummy increments after
full rank, gives
\[
 \Pr[\text{the first }h\text{ columns have rank less than }k]
 \leq \Pr[\operatorname{Bin}(h,1/2)<k]
 \leq \exp(-h/16).
\]
Here $h\geq4k$, so the last estimate follows from Chernoff's inequality.

Suppose the first $h$ columns have rational rank $k$. Select $k$
independent columns by a deterministic rule and denote their nonzero
integer determinant by $b$. Hadamard's inequality gives
$|b|\leq k^{k/2}$, so the number of distinct primes dividing $b$ is at
most $(k/2)\log_2\max(2,k)$. Conditional on the first $h$ columns,
this list of primes is fixed. For every prime not on the list, the
first $h$ columns already have full rank modulo that prime. For each
prime on the list, the remaining $m-h$ columns, independently of the
first block, have full rank modulo that prime except with probability
at most $\exp(-h/16)$, by the preceding argument.

A union bound therefore proves that $B$ has full rank modulo every
prime except with the stated probability. An integer matrix has image
$\mathbb Z^k$ precisely when it has full rank modulo every prime:
this follows immediately from its Smith normal form. This proves the
claim.
\end{proof}

\subsection{Degree classes and isolated extreme degrees}

All logarithms below, unless a base is indicated, are natural.

\begin{lemma}\label{largeprime:binomial}
Fix $a>0$. Suppose $N=n+O(\sqrt n)$ and
\[
 u=\frac N2-a\sqrt{n\log n}+O(\sqrt n).
\]
Then
\[
 \Pr[\operatorname{Bin}(N,1/2)=\lfloor u\rfloor]
     =n^{-1/2-2a^2+o(1)},\qquad
 \Pr[\operatorname{Bin}(N,1/2)\leq\lfloor u\rfloor]
     =n^{-2a^2+o(1)}.
\]
The error is uniform when the constants in the two $O(\sqrt n)$
terms are bounded. The analogous upper-tail estimates hold by
symmetry.
\end{lemma}

\begin{proof}
Put $x=N/2-\lfloor u\rfloor$. Stirling's formula, followed by expansion
of the binomial entropy at $1/2$, gives, in the present range,
\[
 \log\Pr[\operatorname{Bin}(N,1/2)=\lfloor u\rfloor]
 =-\frac12\log N-\frac{2x^2}{N}
    +O\left(1+\frac{x^4}{N^3}+\frac{|x|}{N}\right).
\]
Here $2x^2/N=2a^2\log n+o(\log n)$ uniformly, proving the local estimate.
Hoeffding's inequality gives the required tail upper bound. For the
lower bound, sum the local probabilities at the
$\lfloor\sqrt n/\log n\rfloor$ consecutive integers ending at
$\lfloor u\rfloor$. Each summand has the stated local exponent,
uniformly, and the number of summands is $n^{1/2-o(1)}$.
\end{proof}

\begin{lemma}\label{largeprime:degree-count}
For an integer $u$, let
\[
 Z(u)=|\{v:\deg_G(v)\leq u\}|,\qquad
 \pi=\Pr[\operatorname{Bin}(n-1,1/2)\leq u],\qquad
 b=\Pr[\operatorname{Bin}(n-2,1/2)=u].
\]
Then
\[
 \frac{\operatorname{Var} Z(u)}{(\mathbb E Z(u))^2}
 \leq\frac1{n\pi}+\frac{b^2}{4\pi^2}
\]
whenever $\pi>0$. An analogous bound holds for upper-tail vertex counts,
with the corresponding boundary point probability.
\end{lemma}

\begin{proof}
For distinct vertices $i,j$, write
\[
 \deg_G(i)=A_{ij}+U_i,\qquad \deg_G(j)=A_{ij}+U_j.
\]
The variables $U_i,U_j$ are independent
$\operatorname{Bin}(n-2,1/2)$ variables and are independent of $A_{ij}$.
Conditional on $A_{ij}$, the two degree events are independent. The
two possible conditional probabilities differ by exactly $b$.
Their covariance is therefore $b^2/4$. Since
$\mathbb E Z(u)=n\pi$, summing variances and covariances proves the bound.
The same argument applies to upper tails.
\end{proof}

For the remainder of the proof, set
\begin{equation}\label{largeprime:parameters}
 \begin{aligned}
 k&=\lfloor n^{0.02}\rfloor,&
 K&=\lceil10\sqrt n\rceil,&
 s&=\lfloor n^{0.505}\rfloor,\\
 r&=\lfloor n^{0.14}\rfloor,&
 g&=\lfloor n^{0.03}\rfloor,&
 t&=\left\lfloor\frac{n-1}{2}
             -\sqrt{0.245\,n\log n}\right\rfloor.
 \end{aligned}
\end{equation}
All subsequent assertions involving these parameters are for
sufficiently large $n$. The two size bounds $k$ and $K$ apply,
respectively, to low-degree classes and to arbitrary degree classes.
The $r$ vertices of largest degree will supply the initial coordinates;
the second round will use $s$ low-degree vertices. The inequalities
\[
 k\log n=o(r),\qquad K\log n=o(s),\qquad k=o(g)
\]
allow union bounds over the row sets and ensure that the initial
vertices remain identifiable after deleting at most $k$ vertices.
The threshold $t$ gives $n^{0.51+o(1)}$ low-degree vertices, more
than the $s=n^{0.505+o(1)}$ needed in the second round.

\begin{lemma}\label{largeprime:degree-properties}
With probability $1-o(1)$, all three assertions hold:
\begin{enumerate}
 \item every degree class has fewer than $K$ vertices;
 \item every degree class of degree at most $t$ has fewer than $k$
       vertices;
 \item at least $2s$ vertices have degree at most $t-K$.
\end{enumerate}
\end{lemma}

\begin{proof}
Fix a set $C$ of $K$ vertices and an integer $d$. Conditional on the
edges internal to $C$, the external degrees of its vertices are
independent $\operatorname{Bin}(n-K,1/2)$ variables. Each required
external degree has probability at most $2/\sqrt n$. Thus
\[
 \Pr[\deg_G(v)=d\text{ for all }v\in C]\leq(2/\sqrt n)^K.
\]
Taking a union bound over the $n$ possible degrees and all $K$-sets
gives
\[
 \Pr[\text{assertion 1 fails}]
 \leq n\binom nK(2/\sqrt n)^K
 \leq n(2e\sqrt n/K)^K
 =\exp(-\Omega(\sqrt n)).
\]

For assertion 2, fix a $k$-set $C$ and again condition on its internal
edges. If the common full degree is $d\leq t$, the required external
degree of each vertex is at most $t$. The external degree has
distribution $\operatorname{Bin}(n-k,1/2)$. Since $k=o(\sqrt n)$,
Lemma~\ref{largeprime:binomial} and monotonicity of the binomial
point probabilities to the left of the mode give the uniform bound
$n^{-0.99+o(1)}$ for every required external degree. The external
degrees of the vertices of $C$ are independent. Consequently,
\[
 \Pr[\text{assertion 2 fails}]
 \leq n\binom nk n^{(-0.99+o(1))k}
 \leq n\left(\frac{e n^{0.01+o(1)}}{k}\right)^k
 =\exp(-(0.01-o(1))k\log n)
 =o(1).
\]

For assertion 3, put $u=t-K$ in Lemma~\ref{largeprime:degree-count}.
Lemma~\ref{largeprime:binomial} gives
$\pi=n^{-0.49+o(1)}$ and $b=n^{-0.99+o(1)}$. Hence
\[
 \mathbb E Z(u)=n^{0.51+o(1)},\qquad
 \frac{\operatorname{Var}Z(u)}{(\mathbb E Z(u))^2}
 \leq n^{-0.51+o(1)}+n^{-1+o(1)}=o(1).
\]
Since $\mathbb E Z(u)/s\to\infty$, Chebyshev's inequality gives
$Z(u)\geq\mathbb E Z(u)/2>2s$ except with probability at most
$n^{-0.51+o(1)}+n^{-1+o(1)}$.
The uniform estimate includes the shift by $K=O(\sqrt n)$.
\end{proof}

\begin{lemma}\label{largeprime:seeds}
With failure probability at most $n^{-0.15+o(1)}$, the $r+1$ largest
degrees are pairwise
separated by more than $g$. Let $S$ consist of the $r$ vertices of
largest degree on this event. Then each vertex of $S$ is a singleton
degree class. Moreover, if $C\cap S=\varnothing$ and $|C|\leq k$,
the $r$ vertices of largest degree in $G[[n]\setminus C]$ are
exactly $S$.
\end{lemma}

\begin{proof}
Set
\[
 H=\left\lceil\frac{n-1}{2}+\sqrt{0.42\,n\log n}\right\rceil.
\]
Lemmas~\ref{largeprime:binomial} and~\ref{largeprime:degree-count}
show that the number of vertices with degree at least $H$ has
expectation $n^{0.16+o(1)}$ and relative variance at most
$n^{-0.16+o(1)}+n^{-1+o(1)}$.
By Chebyshev's inequality it is at least $r+1$, except with probability
$n^{-0.16+o(1)}$.

For a fixed pair of distinct vertices $i,j$, condition on their
mutual edge. The remaining parts of their degrees are independent
$\operatorname{Bin}(n-2,1/2)$ variables. Their upper-tail probabilities
at the relevant threshold are at most $n^{-0.84+o(1)}$, and their
largest point probabilities above that threshold are at most
$n^{-1.34+o(1)}$. Conditional on the first remaining degree, at most
$2g+1$ values of the second remaining degree are within distance $g$.
It follows that
\[
 \Pr[\deg_G(i),\deg_G(j)\geq H,\
          |\deg_G(i)-\deg_G(j)|\leq g]
 \leq (2g+1)n^{-0.84+o(1)}n^{-1.34+o(1)}
 =n^{-2.15+o(1)}.
\]
A union bound over all pairs gives failure probability
$n^{-0.15+o(1)}=o(1)$ for the required separation above $H$.
Together with the preceding count, this proves the first assertion.

The separation includes the $(r+1)$st largest degree, so no vertex
of $S$ shares its degree with a vertex outside $S$. Deleting a
disjoint set $C$ of size at most $k$ lowers every remaining degree
by an integer between $0$ and $k$. Since $g>k$, every vertex of $S$
still has degree greater than every vertex outside $S\cup C$ in
the induced graph. This proves the deletion assertion.
\end{proof}

\begin{remark}
The use of separated extreme degrees goes back to
Babai, Erd\H{o}s, and Selkow~\cite[Theorem~3.6]{BabaiErdosSelkow1980},
who prove separation of the largest
$\lfloor n^{0.15}\rfloor$ degrees on the scale
$\lfloor n^{0.03}\rfloor$. We included the tail and pair estimates
to record the failure probability and the stability under deletion
needed for the integer construction.
\end{remark}

\subsection{Selecting independent adjacency blocks}

The next two lemmas are proved for all small row sets, without
requiring them to be degree classes. This order of quantifiers permits
their later use for the random degree partition.

\begin{lemma}\label{largeprime:first-blocks}
For every fixed nonempty $C\subseteq[n]$ with $|C|\leq k$, let
$S_C$ be the $r$ highest-degree vertices of $G[[n]\setminus C]$,
using vertex labels to break ties. With probability $1-o(1)$,
\[
 A[C,S_C]\mathbb Z^r=\mathbb Z^C
 \quad\text{for every such }C.
\]
\end{lemma}

\begin{proof}
Fix $C$ and condition only on $G[[n]\setminus C]$. Then $S_C$
is fixed, and the entries of $A[C,S_C]$ are independent fair
Bernoulli variables: they concern edges with one endpoint in $C$,
none of which has been exposed. Lemma~\ref{largeprime:rectangular}
applies because $r\geq16k$ for sufficiently large $n$.

There are at most $(k+1)n^k$ eligible sets $C$. The lemma's polynomial
prefactor is absorbed into the following estimate, and a union bound
gives total failure probability
\[
 \exp\bigl(O(k\log n)-\Omega(r)\bigr)=o(1).
\]
\end{proof}

\begin{lemma}\label{largeprime:second-blocks}
For every fixed nonempty $C\subseteq[n]$ with $|C|\leq K$, define
\[
 R_C=\{v\notin C:\deg_{G[[n]\setminus C]}(v)\leq t-K\}.
\]
If $|R_C|\geq s$, let $T_C$ be its first $s$ vertices in label order.
With probability $1-o(1)$, simultaneously for every such $C$ with
$|R_C|\geq s$,
\[
 A[C,T_C]\mathbb Z^s=\mathbb Z^C.
\]
On assertion 3 of Lemma~\ref{largeprime:degree-properties},
$|R_C|\geq s$ for every eligible $C$. Every vertex of $R_C$
has full degree at most $t$.
\end{lemma}

\begin{proof}
Fix $C$ and condition on $G[[n]\setminus C]$. The event $|R_C|\geq s$
and, when it exists, the set $T_C$, are determined by this conditioning.
The entries of $A[C,T_C]$ remain independent fair Bernoulli variables.
Since $s\geq16K$ for sufficiently large $n$,
Lemma~\ref{largeprime:rectangular} and a union bound over at most
$(K+1)n^K$ sets give total failure probability
\[
 \exp\bigl(O(K\log n)-\Omega(s)\bigr)=o(1).
\]
Indeed $s/(\sqrt n\log n)\longrightarrow\infty$.

If a vertex outside $C$ has full degree at most $t-K$, deleting
$C$ cannot increase its degree, so that vertex belongs to $R_C$.
On the stated degree-count event,
\[
 |R_C|\geq|\{v:\deg_G(v)\leq t-K\}|-|C|
          \geq2s-K\geq s.
\]
Conversely, for $v\in R_C$, restoring the edges to $C$ increases its
degree by at most $|C|\leq K$, and therefore $\deg_G(v)\leq t$.
\end{proof}

After averaging over the exposed induced graph, both bounds are
unconditional. We may therefore intersect these uniform events with
the degree and seed events and then specialize $C$ to an actual
degree class.

\subsection{The integer bootstrap}

\begin{proof}[Proof of Theorem~\ref{largeprime:main}]
Intersect the events in
Lemmas~\ref{largeprime:degree-properties}--\ref{largeprime:second-blocks}.
The degree-count event fails with probability at most
$n^{-0.51+o(1)}+n^{-1+o(1)}$; the seed event fails with probability at
most $n^{-0.15+o(1)}$; and the other failure probabilities decrease
faster than any fixed inverse power of $n$.
Thus the intersection has failure probability at most
$n^{-0.15+o(1)}$. We prove the claimed
integer inclusion for every graph in it.

For $d\in\mathcal D_G$, let $P_d$ be the coordinate projector onto the
vertices of degree $d$. Lagrange interpolation gives
\[
 P_d=\prod_{\substack{h\in\mathcal D_G\\h\neq d}}
                      \frac{D-hI}{d-h}.
\]
The denominator product divides $\Delta_G$ up to sign, so
\begin{equation}\label{largeprime:integer-projector}
 \Delta_G P_d\in\mathbb Z[D].
\end{equation}
Thus this operator preserves $\Lambda_G$. For a degree class $C$,
we also write $P_C$ for its projector.

Each $v\in S$ is a singleton degree class, so
\[
 \Delta_G e_v=(\Delta_G P_{\deg_G(v)})e\in\Lambda_G.
\]
Now let $C$ be an actual degree class of size at most $k$.
If it meets $S$, it is a singleton already covered. Otherwise,
Lemma~\ref{largeprime:seeds} gives $S_C=S$, and
Lemma~\ref{largeprime:first-blocks} gives
$A[C,S]\mathbb Z^r=\mathbb Z^C$. For every $v\in S$,
\[
 (\Delta_G P_C)A(\Delta_G e_v)
       =\Delta_G^2 P_CAe_v\in\Lambda_G.
\]
Integer combinations of these columns therefore give
$\Delta_G^2 e_u\in\Lambda_G$ for every $u\in C$. For classes meeting
$S$, the same assertion follows by multiplying the already
available vector by $\Delta_G$.
By assertion 2 of Lemma~\ref{largeprime:degree-properties}, we have
proved
\begin{equation}\label{largeprime:first-round}
 \Delta_G^2 e_u\in\Lambda_G
 \quad\text{whenever }\deg_G(u)\leq t.
\end{equation}

Finally let $C$ be any actual degree class. Its size is less than $K$.
Lemma~\ref{largeprime:second-blocks} supplies a set $T_C$ of $s$
vertices, each of full degree at most $t$, such that
$A[C,T_C]\mathbb Z^s=\mathbb Z^C$. By
\eqref{largeprime:first-round} and \eqref{largeprime:integer-projector},
for every $v\in T_C$,
\[
 (\Delta_G P_C)A(\Delta_G^2 e_v)
       =\Delta_G^3P_CAe_v\in\Lambda_G.
\]
Integer combinations give $\Delta_G^3e_u\in\Lambda_G$ for all $u\in C$.
The degree classes partition $[n]$, proving
$\Delta_G^3\mathbb Z^n\subseteq\Lambda_G$.

Each operator $\Delta_GP_d$ is an integer polynomial in $D$ of
degree at most $q-1$, so we can also track the words used.
The initial seed vectors use no $A$, the first round adds one $A$,
and the second adds one more. Thus the entire construction uses
integer combinations of
\[
 D^i e,\qquad D^iAD^j e,\qquad D^iAD^jAD^\ell e,
 \qquad 0\leq i,j,\ell\leq q-1.
\]
It therefore gives the stronger inclusion
$\Delta_G^3\mathbb Z^n\subseteq\Lambda_G^{[2]}$ asserted in the
theorem. The argument also covers initially seeded vertices used in
the first round by merely multiplying their available vectors by
the integer $\Delta_G$.

This inclusion implies that the cokernel is finite and is annihilated
by $\Delta_G^3$. For any prime $p\nmid\Delta_G$, reducing the inclusion
modulo $p$ gives full span over $\mathbb F_p$. Every prime divisor
of $\Delta_G$ divides a nonzero degree difference, and is therefore
at most $w_G$. These conclusions hold on the same graph event.
\end{proof}

\begin{corollary}\label{largeprime:width}
With probability $1-o(1)$, the AD word module of $G(n,1/2)$ is full
over $\mathbb F_p$ simultaneously for every prime
$p>2\sqrt{n\log n}$.
\end{corollary}

\begin{proof}
Hoeffding's inequality and a union bound imply
\[
 \Pr\left[\max_v\left|\deg_G(v)-\frac{n-1}{2}\right|
                  >\sqrt{n\log n}\right]\leq\frac2n.
\]
Hence $w_G\leq2\sqrt{n\log n}$ with probability $1-o(1)$, and
Theorem~\ref{largeprime:main} applies.
\end{proof}

\begin{remark}\label{largeprime:remaining}
Theorem~\ref{largeprime:main} reduces integral saturation to finitely
many primes: every possible prime divisor of the defect is at most
$w_G\leq n-1$. We next treat the prime two, each fixed odd prime,
and a uniform interval of growing primes. These estimates complete
the proof of Theorem~\ref{integral:random}.
\end{remark}

%% file: article_random_integral.tex
\subsection{Completion of the proof of Theorem~\ref{integral:random}}
\label{integral:section}

Fix the absolute cutoff \(P_0\) in Theorem~\ref{growing:uniform}.
Theorem~\ref{largeprime:main} gives
\[
 \Delta_G^3\mathbb Z^n\subseteq\Lambda_G
\]
except with probability at most \(n^{-3/20+o(1)}\).
Every prime divisor of \(\Delta_G\) divides a nonzero difference
between two degrees in \(\{0,\ldots,n-1\}\), and is therefore at most
\(n-1\). On this event, \(M_p(G)=\mathbb F_p^n\) simultaneously
for every prime \(p>n\).

It remains to treat primes at most \(n\). At the prime two,
Theorem~\ref{modtwo:random} gives failure probability
\(\exp(-\Omega(n))\). For the finitely many odd primes below \(P_0\),
Theorem~\ref{odd:fixed-irreducible} and a finite union bound give
total failure probability \(O(n^{-1})\). Finally,
Theorem~\ref{growing:uniform} covers all primes \(P_0\leq p\leq n\)
with failure probability \(\exp(-\Omega(\sqrt n))\).
The cutoff was fixed before taking \(n\to\infty\), so the constants
in the fixed-prime estimates remain fixed.

On the intersection of these events, \(\Lambda_G\) has full image
over every prime field. Proposition~\ref{prop:localcriterion}
therefore gives \(\Lambda_G=\mathbb Z^n\), and
\[
 \begin{aligned}
 \Pr[\Lambda_G\neq\mathbb Z^n]
 &\leq n^{-3/20+o(1)}+\exp(-\Omega(n))
           +O(n^{-1})+\exp(-\Omega(\sqrt n))\\
 &=n^{-3/20+o(1)}.
 \end{aligned}
\]
This proves Theorem~\ref{integral:random}. The exponent records
the degree thresholds used in Section~\ref{largeprime:section};
we have not attempted to optimize it.

%% file: article_caterpillar_reconstruction.tex
\section{Reconstruction from small-diameter caterpillars}
\label{sec:caterpillar-reconstruction}

The principal-moment interpretation of Lu and Li~\cite[Section~7]{LuLi2026}
gives a combinatorial form of our rigidity result. For nonnegative integers
\(a_0,\ldots,a_k\), let \(C(a_0,\ldots,a_k)\) be the tree obtained by
attaching \(a_i\) leaves at vertex \(i\) of a \(k\)-edge path. Then
\begin{equation}\label{eq:caterpillar-dictionary}
 \operatorname{hom}(C(a_0,\ldots,a_k),G)
     =e^{\mathsf T}D^{a_0}AD^{a_1}\cdots AD^{a_k}e.
\end{equation}
Here homomorphisms need not be injective, and the test trees carry no
labels or weights. Every AD word has this form, allowing zero exponents.
A word of length \(L\) with \(k\) occurrences of \(A\) gives a
caterpillar with \(L+1\) vertices and diameter at most \(k+2\).
The identity and the finite Gram-matrix reconstruction method are
established ingredients. The extra conclusion here is that a fixed
family of test trees of bounded diameter and sublinear order suffices
for almost every graph.

We first record the rigidity property of the integer construction.
It explains why this reconstruction consequence does not require the
local-prime arguments.

\begin{remark}[Rigidity uses less than integral saturation]
\label{rem:implications-shorter-rigidity}
Let \(\mathcal E_n\) be the intersection of the events in
Lemmas~\ref{largeprime:degree-properties}--\ref{largeprime:second-blocks}.
The proof of Theorem~\ref{largeprime:main} gives
\(\Pr(\mathcal E_n^c)\leq n^{-3/20+o(1)}\).
Every real orthogonal pointed AD intertwiner from a graph in
\(\mathcal E_n\) to any graph \(H\) is a permutation matrix.

Indeed, such an intertwiner \(Q\) maps each degree block of \(G\)
to the corresponding block of \(H\). A singleton degree vertex \(s\)
therefore satisfies \(Qe_s=e_{s'}\): orthogonality gives a sign,
and \(Qe=e\) makes it positive. Suppose a set \(S\) of vertices
has already been recovered in this sense, and
\(A_G[C,S]\mathbb Z^S=\mathbb Z^C\) for a degree block \(C\).
Each \(e_u\), \(u\in C\), is an integer combination of the
vectors \(P_CA_Ge_s\). Their images are
\(P_{C'}A_He_{s'}\), which are integer vectors. Hence \(Qe_u\)
is integral; orthogonality and pointedness make it a positive
coordinate unit vector. The two rounds of the construction recover
every vertex.

On the same event, the inclusion
\(\Delta_G^3\mathbb Z^n\subseteq\Lambda_G^{[2]}\) gives AD-fullness,
so Theorem~\ref{thm:finiteGNS} supplies \(Q\) when the principal
moments agree. The rigidity argument uses the actual cut conditions
defining \(\mathcal E_n\), not merely that inclusion. It also does
not prove integral saturation: the degree projectors need not lie
in \(\mathbb Z[D]\). The local-prime arguments remove these
remaining interpolation denominators.
\end{remark}

\begin{corollary}[A fixed family of small-diameter tests]
\label{cor:random-caterpillars}
Put
\[
 m_n=\left\lceil12\sqrt{n\log n}+6\right\rceil,
 \qquad
 \mathcal F_n=\{T:T\text{ is a caterpillar},\
           \operatorname{diam}T\leq7,\ |V(T)|\leq m_n\},
\]
including \(K_1\). If \(G\sim G(n,1/2)\), then with probability
at least \(1-n^{-3/20+o(1)}\), every finite simple graph \(H\) satisfies
\[
 \bigl[\operatorname{hom}(T,G)=\operatorname{hom}(T,H)
                   \text{ for all }T\in\mathcal F_n\bigr]
       \quad\Longrightarrow\quad G\cong H.
\]
The family \(\mathcal F_n\) has \(O(n^3(\log n)^3)\) isomorphism types.
\end{corollary}

\begin{proof}
Fix \(G\in\mathcal E_n\), and let \(q\) be its number of distinct
degrees. The construction in Theorem~\ref{largeprime:main} shows that
the vectors obtained from the words
\[
 D^i,\qquad D^iAD^j,\qquad D^iAD^jAD^\ell,
       \qquad 0\leq i,j,\ell<q,
\]
span \(\mathbb R^n\). Choose a basis among them, with the empty
word first, as the columns of a matrix \(B_G\). Evaluate the same
words on \(H\) to obtain \(B_H\).
Every chosen word has length at most \(3q-1\) and at most two
occurrences of \(A\).

Suppose first that \(G,H\) have equal principal moments for all
words of length at most \(6q-1\) containing at most five \(A\)'s.
The empty word gives \(|V(H)|=n\). As in the finite reconstruction
argument of Theorem~\ref{thm:finiteGNS}, the stated moments give
\[
 B_G^{\mathsf T}B_G=B_H^{\mathsf T}B_H,\qquad
 B_G^{\mathsf T}X_GB_G=B_H^{\mathsf T}X_HB_H
       \quad (X=A,D).
\]
The Gram matrix is positive definite, so \(B_H\) is invertible.
These identities show that \(Q=B_HB_G^{-1}\) is orthogonal
and intertwines both \(A\) and \(D\). The first columns give \(Qe=e\).
Remark~\ref{rem:implications-shorter-rigidity} makes \(Q\) a
permutation, proving \(G\cong H\). No good-event assumption on
\(H\) is needed.

By \eqref{eq:caterpillar-dictionary}, all the required moments
are counts from caterpillars with at most \(6q\) vertices and
diameter at most seven. Hoeffding's inequality and a union bound give
\[
 \Pr\!\left[\max_v\left|\deg_G(v)-\frac{n-1}{2}\right|
                    >\sqrt{n\log n}\right]\leq\frac2n.
\]
On the complementary event, \(q\leq2\sqrt{n\log n}+1\), so
\(6q\leq m_n\). Combining the two exceptional probabilities
proves the claimed bound.

Every caterpillar of diameter at most seven admits a spine with
at most six vertices. For a spine with \(s\) vertices, the number
of nonnegative leaf-multiplicity sequences giving at most \(m_n\)
vertices is \(\binom{m_n}{s}\). Thus the number of isomorphism
types is at most \(\sum_{s=1}^6\binom{m_n}{s}=O(m_n^6)\), as required.
\end{proof}

The next result shows that reducing the test diameter is a substantive
restriction, even when arbitrarily large test trees are allowed.

\begin{proposition}[Diameter four does not suffice]
\label{prop:caterpillar-four}
The probability that \(G\sim G(n,1/2)\) is determined among finite
simple graphs by its homomorphism counts from all caterpillars of
diameter at most four tends to zero.
\end{proposition}

\begin{proof}
Let \(I_n\) be the integer interval centred at \((n-1)/2\) with
radius \(\sqrt{n\log n}\), intersected with \(\{0,\ldots,n-1\}\).
Except with probability at most \(2/n\), all degrees belong to
\(I_n\). Put \(b=|I_n|=O(\sqrt{n\log n})\).
For \(k=0,1,2\) and \(d_0,\ldots,d_k\in I_n\), record
\(N_k(d_0,\ldots,d_k)\), the number of ordered walks
\((v_0,\ldots,v_k)\) with \(\deg(v_i)=d_i\).
There are \(b+b^2+b^3\) such entries, each an integer between
zero and \(n^3\).

Every caterpillar of diameter at most four admits a spine of at
most two edges. Its homomorphism count is therefore a sum of the
entries of the corresponding \(N_k\), weighted by products of
degree powers. Thus the entire infinite family of counts is
determined by these tensors. Among graphs with degrees in \(I_n\),
the number of possible count profiles is at most
\[
 (n^3+1)^{b+b^2+b^3}
   =2^{O(n^{3/2}(\log n)^{5/2})}=2^{o(n^2)}.
\]
A profile that determines a graph up to isomorphism can occur on
at most \(n!\) labelled graphs. Since \(G(n,1/2)\) is uniform
on the \(2^{\binom n2}\) labelled graphs, the desired probability
is at most
\[
 \frac2n+n!\,2^{o(n^2)-\binom n2}=o(1).
\]
\end{proof}

%% file: article_random_mod2.tex
\section{Saturation at the prime two}
\label{modtwo:section}

In this section we work over $\mathbb F_2$ and write
\[
 A=A_G\bmod 2,\qquad
 e=(1,\ldots,1)^{\mathsf T},\qquad
 d=Ae,\qquad
 D=D_G\bmod 2=\operatorname{diag}(d).
\]
Thus the degree-parity diagonal is determined by $A$. The principal
module is
\[
 M_2(G)=\mathbb F_2\langle A,D\rangle e.
\]
All orthogonal complements below refer to the standard
nondegenerate bilinear form $x\cdot y=x^{\mathsf T}y$. For a subspace $W$,
its radical is $W\cap W^\perp$.

\begin{theorem}\label{modtwo:random}
Let $G\sim G(n,1/2)$. For every $n\geq 48$,
\begin{equation}\label{modtwo:failure-bound}
 \Pr\bigl[M_2(G)\neq\mathbb F_2^n\bigr]
 \leq
 4\exp(-n/18)+1024(n+1)^4\,2^{-n/48}.
\end{equation}
Consequently $M_2(G)=\mathbb F_2^n$ with probability
$1-\exp(-\Omega(n))$.
\end{theorem}

We prove the theorem by conditioning on the degree-parity vector.
Failure of cyclicity then gives an invariant subspace orthogonal to
$e$, split between the two parity classes. Counting these subspaces
by dimension and radical makes it possible to compare their number
with the exact cost of invariance in the conditional matrix space.
The numerical constants are not optimized.

\subsection{The degree-parity fibre and invariant obstructions}

Let $\operatorname{Alt}_n(\mathbb F_2)$ be the space of symmetric matrices
with zero diagonal. These are precisely the matrices of alternating
bilinear forms in characteristic two. The adjacency matrix of
$G(n,1/2)$ is uniform on this space.

\begin{lemma}\label{modtwo:degree-fibre}
The vector $d=Ae$ is uniform on $e^\perp$. Conditional on any
$d\in e^\perp$, the matrix $A$ is uniform on the nonempty affine space
\[
 \mathcal F_d=\{A\in\operatorname{Alt}_n(\mathbb F_2):Ae=d\}.
\]
For $i\in\{0,1\}$, put $V_i=\{v:d_v=i\}$ and $n_i=|V_i|$. Then
\begin{equation}\label{modtwo:degree-balance}
 \Pr\bigl[\min(n_0,n_1)<n/3\bigr]\leq 4\exp(-n/18).
\end{equation}
\end{lemma}

\begin{proof}
For every alternating $A$, $e^{\mathsf T}Ae=0$, so $Ae\in e^\perp$.
The image of the edge matrix
$E_{uv}+E_{vu}$ under $A\mapsto Ae$ is $e_u+e_v$.
For $n\geq2$ these vectors span $e^\perp$; for $n=1$ both the
domain and the target image are zero. Thus the linear map is surjective
onto $e^\perp$, and all its fibres have the same size.

The uniform distribution on $e^\perp$ is the distribution of $n$
independent fair bits conditioned on having even sum. The conditioning
event has probability $1/2$ for every $n\geq1$. Hence $n_1$ has the
distribution of a $\operatorname{Bin}(n,1/2)$ variable conditioned to be
even. The usual Hoeffding bound
\[
 \Pr\bigl[|\operatorname{Bin}(n,1/2)-n/2|>n/6\bigr]
 \leq 2\exp(-n/18)
\]
and division by the conditioning probability give
\eqref{modtwo:degree-balance}.
\end{proof}

Fix for now a vector $d\in e^\perp$ such that
$n_0,n_1\geq n/3$. In particular both blocks are nonempty.
Let $e_i$ denote the all-ones vector on the coordinate space
$\mathbb F_2^{V_i}$, and set
\[
 H_i=\{x\in\mathbb F_2^{V_i}:x\cdot e_i=0\}.
\]
When necessary, $e_i$ is also viewed as a vector in $\mathbb F_2^n$ by
extension by zero. With this convention $e=e_0+e_1$ and $d=e_1$.

\begin{lemma}\label{modtwo:invisible}
The equality $M_2(G)\neq\mathbb F_2^n$ holds if and only if there is a
nonzero subspace $W\subseteq e^\perp$ with
\[
 AW\subseteq W,\qquad DW\subseteq W.
\]
Every such $D$-invariant subspace has a unique decomposition
\[
 W=W_0\oplus W_1,\qquad W_i\subseteq H_i.
\]
Write
\[
 t_i=\dim W_i,\qquad
 r_i=\dim(W_i\cap W_i^\perp),\qquad
 \varepsilon_i=\mathbf 1_{\{e_i\in W_i\}},
\]
where the perpendicular in $r_i$ is taken in
$\mathbb F_2^{V_i}$. Then
\[
 t=\dim W=t_0+t_1,\qquad
 r=\dim(W\cap W^\perp)=r_0+r_1,
\]
\begin{equation}\label{modtwo:parameter-ranges}
 1\leq t\leq n-2,\qquad
 0\leq r_i\leq \min(t_i,n_i-t_i).
\end{equation}
Moreover, $e\in W$ if and only if $\varepsilon_0=\varepsilon_1=1$.
\end{lemma}

\begin{proof}
The space $M_2(G)$ is invariant under the symmetric matrices $A,D$.
Therefore its orthogonal complement is invariant under both matrices,
is contained in $e^\perp$, and is nonzero exactly when $M_2(G)$ is
proper. Conversely, if such a nonzero $W$ exists, its orthogonal
complement is invariant under $A,D$ and contains $e$. It therefore
contains $M_2(G)$, which is consequently proper.

The matrices $D$ and $I-D$ are the coordinate projections onto $V_1$
and $V_0$, respectively. Their invariance gives the asserted direct
sum. Each $W_i$, embedded in the full coordinate space, is itself
contained in $W\subseteq e^\perp$, so $W_i\subseteq H_i$.
The two coordinate spaces are orthogonal and nondegenerate; hence
radicals add, and $r_i\leq\min(t_i,n_i-t_i)$.
Since $\dim H_i=n_i-1$ and both blocks are nonempty, $t\leq n-2$.
Finally, the uniqueness of the block decomposition gives the claim
about $e\in W$.
\end{proof}

\subsection{Counting subspaces with a prescribed radical}

On an even-dimensional coordinate block, the dot product restricted
to the even-weight hyperplane has a one-dimensional radical. We
retain this radical in the count; its contribution is needed for
subspaces close to full dimension.

Put
\[
 \gamma=\prod_{j=1}^{\infty}(1-2^{-j})^{-1}<4.
\]
Indeed, $\prod_{j=1}^3(1-2^{-j})=21/64$, and the tail product is
at least $1-\sum_{j\geq4}2^{-j}=7/8$.
The Gaussian binomial coefficient satisfies
\begin{equation}\label{modtwo:gaussian}
 \genfrac{[}{]}{0pt}{}{N}{a}_2\leq \gamma\,2^{a(N-a)}
 \qquad (0\leq a\leq N).
\end{equation}

\begin{lemma}\label{modtwo:symplectic-count}
Let $S$ be a nondegenerate symplectic space of dimension $h$ over
$\mathbb F_2$. The number $\mathcal N(h;a,b)$ of subspaces of dimension $a$
and radical dimension $b$ is at most
\begin{equation}\label{modtwo:symplectic-upper}
 \mathcal N(h;a,b)\leq
 \gamma^2\,2^{a(h-a)-b(b-1)/2}.
\end{equation}
Impossible parameter choices contribute zero. More explicitly, for
feasible parameters,
\begin{equation}\label{modtwo:symplectic-exact}
 \mathcal N(h;a,b)=
 I(h,b)\,
 \frac{|\operatorname{Sp}(h-2b,2)|}
 {|\operatorname{Sp}(a-b,2)|\,|\operatorname{Sp}(h-a-b,2)|},
\end{equation}
where
\begin{equation}\label{modtwo:isotropic-count}
 I(h,b)=
 \prod_{j=0}^{b-1}\frac{2^{h-j}-2^j}{2^b-2^j},
 \qquad
 |\operatorname{Sp}(2s,2)|
 =2^{s^2}\prod_{j=1}^{s}(2^{2j}-1).
\end{equation}
Empty products and the order of $\operatorname{Sp}(0,2)$ are one.
\end{lemma}

\begin{proof}
A radical $R$ of dimension $b$ is a totally isotropic subspace.
After an ordered isotropic independent $j$-tuple has been chosen,
its perpendicular space has $2^{h-j}$ vectors, of which $2^j$ lie
in its span. Counting ordered isotropic bases and dividing by the
number of ordered bases of $\mathbb F_2^b$ gives
\eqref{modtwo:isotropic-count}. In particular,
\begin{equation}\label{modtwo:isotropic-upper}
 I(h,b)\leq
 \gamma\,2^{b(h-b)-b(b-1)/2}.
\end{equation}

For fixed $R$, any subspace $U$ with radical exactly $R$ satisfies
$R\subseteq U\subseteq R^\perp$. Its image $U/R$ is a nondegenerate
subspace of dimension $a-b$ in the nondegenerate symplectic space
$R^\perp/R$ of dimension $h-2b$. Conversely every such image has
a unique inverse image with radical $R$. The symplectic group is
transitive on nondegenerate subspaces of a fixed dimension: a
symplectic basis of the subspace extends by a symplectic basis of
its perpendicular complement. The stabilizer is the product of
the symplectic groups on those two summands. This proves
\eqref{modtwo:symplectic-exact}; the group order in
\eqref{modtwo:isotropic-count} follows by counting symplectic bases.

For the upper bound, it is enough to discard the nondegeneracy
condition on $U/R$ and use \eqref{modtwo:gaussian}. Combining it
with \eqref{modtwo:isotropic-upper} gives an exponent
\[
 b(h-b)-\frac{b(b-1)}2
 +(a-b)(h-a-b)
 =a(h-a)-\frac{b(b-1)}2,
\]
which proves \eqref{modtwo:symplectic-upper}.
\end{proof}

\begin{lemma}\label{modtwo:even-hyperplane-count}
Let $V=\mathbb F_2^N$, let $e_N$ be its all-ones vector, and let
$H=e_N^\perp$, with $N\geq1$.
For $W\subseteq H$, write
\[
 t=\dim W,\qquad r=\dim(W\cap W^\perp),\qquad
 \varepsilon=\mathbf 1_{\{e_N\in W\}}.
\]
The number of such $W$ for fixed $(t,r,\varepsilon)$ is at most
\begin{equation}\label{modtwo:block-count}
 \gamma^2\,2^{E(N,t,r,\varepsilon)},\qquad
 E(N,t,r,\varepsilon)
 =t(N-1-t)-\frac{r(r-1)}2
       -\varepsilon(N-t-r).
\end{equation}
The exact counts, using the notation of
Lemma~\ref{modtwo:symplectic-count}, are
\[
 \begin{cases}
  \mathcal N(N-1;t,r),& N\text{ odd},\ \varepsilon=0,\\
  2^t \mathcal N(N-2;t,r),& N\text{ even},\ \varepsilon=0,\\
  \mathcal N(N-2;t-1,r-1),& N\text{ even},\ \varepsilon=1,\\
  0,&N\text{ odd},\ \varepsilon=1.
 \end{cases}
\]
\end{lemma}

\begin{proof}
For $v\in H$,
$v\cdot v=\sum_jv_j=v\cdot e_N=0$, so the restriction of the
dot product to $H$ is alternating. Its radical is
$H\cap\langle e_N\rangle$.

If $N$ is odd, then $e_N\notin H$, and $H$ is nondegenerate
symplectic of dimension $N-1$. The first and fourth cases follow.

Suppose $N$ is even. Then the radical of $H$ is
$L=\langle e_N\rangle$, and $\overline H=H/L$ is nondegenerate
symplectic of dimension $N-2$.
If $\varepsilon=0$, the projection of $W$ into $\overline H$ is
injective, preserving the restricted bilinear form. Its image has
dimension $t$ and radical dimension $r$. For each fixed image,
choose any linear splitting of its inverse image in $H$.
The subspaces projecting isomorphically onto it are exactly the
graphs of linear maps from a $t$-dimensional space to $L$.
There are $2^t$ such maps. This proves the second exact formula.
Its exponent in \eqref{modtwo:symplectic-upper}, after adding $t$,
is $t(N-1-t)-r(r-1)/2$.

If $\varepsilon=1$, then $L\subseteq W$. Its image has dimension
$t-1$ and radical dimension $r-1$, and $W$ is its unique full
inverse image. The resulting exponent is
\[
 (t-1)(N-t-1)-\frac{(r-1)(r-2)}2
 =t(N-1-t)-\frac{r(r-1)}2-(N-t-r).
\]
This proves the last remaining case and the bound.
\end{proof}

Apply Lemma~\ref{modtwo:even-hyperplane-count} independently in
the two degree-parity blocks. For a fixed choice of the six
parameters $t_0,t_1,r_0,r_1,\varepsilon_0,\varepsilon_1$,
the number of candidate subspaces $W=W_0\oplus W_1$ is at most
\begin{equation}\label{modtwo:two-block-count}
 \gamma^4\,2^{E_0+E_1},\qquad
 E_i=t_i(n_i-1-t_i)-\frac{r_i(r_i-1)}2
                  -\varepsilon_i(n_i-t_i-r_i).
\end{equation}

\subsection{Invariance inside a fixed degree-parity fibre}

\begin{lemma}\label{modtwo:conditional-codimension}
Let $V=\mathbb F_2^n$, let $e\neq0$, and fix $d\in e^\perp$.
Let $W\subseteq e^\perp$ have dimension $t$ and radical dimension $r$.
For a matrix $A$ chosen uniformly from
$\{A\in\operatorname{Alt}_n(\mathbb F_2):Ae=d\}$, put
\[
 c=t(n-t)-\frac{r(r+1)}2.
\]
If $e\notin W$, the conditions $AW\subseteq W$ and $Ae=d$ are
compatible if and only if $d\in W^\perp$. When compatible,
\[
 \Pr[AW\subseteq W\mid Ae=d]=2^{-(c-t)}.
 \]
If $e\in W$, they are compatible if and only if
$d\in W\cap W^\perp$. When compatible,
\[
 \Pr[AW\subseteq W\mid Ae=d]=2^{-(c-n+r+1)}.
\]
In an incompatible case the probability is zero.
\end{lemma}

\begin{proof}
The dual space of alternating forms identifies with $\bigwedge^2V$:
the bivector $u\wedge v$ evaluates the form $A$ as $u^{\mathsf T}Av$.
This is valid also in characteristic two, where $u\wedge u=0$.
Since the dot product is nondegenerate, $AW\subseteq W$ is
equivalent to
\[
 v^{\mathsf T}Aw=0\qquad(w\in W,\ v\in W^\perp).
\]
Thus its homogeneous constraint space is
$\mathcal U=W\wedge W^\perp\subseteq\bigwedge^2V$.

Write $R=W\cap W^\perp$, and choose a direct sum
\[
 V=R\oplus P\oplus Q\oplus Z,\qquad
 W=R\oplus P,\qquad W^\perp=R\oplus Q.
\]
Here $\dim P=t-r$, $\dim Q=n-t-r$, and $\dim Z=r$.
The independent wedge types in $\mathcal U$ are
\[
 \bigwedge\nolimits^2R,\quad R\wedge P,\quad
 R\wedge Q,\quad P\wedge Q.
\]
Their dimensions give
\[
 \dim\mathcal U
 =\binom r2+r(t-r)+r(n-t-r)+(t-r)(n-t-r)
 =c.
\]

The degree condition has dual constraint space
$\mathcal E_e=e\wedge V$, of dimension $n-1$. It prescribes the
value $d\cdot v$ on $e\wedge v$; the assumption $d\cdot e=0$
makes this well-defined.
Because $W\subseteq e^\perp$, we have $e\in W^\perp$.
If $e\notin W$, choose $Q$ to contain $e$. An adapted wedge basis
then gives
\[
 \mathcal U\cap\mathcal E_e=e\wedge W,
 \qquad \dim(\mathcal U\cap\mathcal E_e)=t.
\]
Indeed the wedges of $e$ with $R$ and $P$ are present, whereas
wedges with the remaining basis vectors in $Q$ or $Z$ are absent.
If $e\in W$, then $e\in R$, and the same basis argument yields
\[
 \mathcal U\cap\mathcal E_e=e\wedge(W+W^\perp),
 \qquad
 \dim(\mathcal U\cap\mathcal E_e)=n-r-1.
\]

The affine degree conditions and the zero conditions on
$\mathcal U$ are compatible precisely when their prescribed
values agree on this intersection. In the first case this is
$d\cdot W=0$. In the second case it is
$d\cdot(W+W^\perp)=0$, equivalently
$d\in(W+W^\perp)^\perp=W\cap W^\perp$.
When compatible, the additional independent constraints on the
degree fibre number $\dim\mathcal U-
\dim(\mathcal U\cap\mathcal E_e)$. A uniform point of an affine
space over $\mathbb F_2$ satisfies that many independent affine
constraints with probability $2^{-\text{codimension}}$.
This proves the two formulas.
\end{proof}

For the graded subspaces of Lemma~\ref{modtwo:invisible}, compatibility
is automatic. If $e\notin W$, then $d=e_1$ is orthogonal to $W$
because $W_1\subseteq H_1$. If $e\in W$, both $\varepsilon_i$ are
one, so $e_1\in W$ as well, and $d=e_1\in W\cap W^\perp$.
In particular, this second case can occur only when both $n_i$
are even.

\subsection{Balancing entropy and codimension}

\begin{proof}[Proof of Theorem~\ref{modtwo:random}]
Fix a balanced degree-parity vector $d$.
For any feasible nonzero candidate $W$, set
\[
 E=E_0+E_1,\qquad
 X=t_0(n_1-t_1)+t_1(n_0-t_0).
\]
Let $c_{\mathrm{cond}}$ be its conditional invariance codimension
from Lemma~\ref{modtwo:conditional-codimension}.
Substitution of \eqref{modtwo:two-block-count} gives, if $e\notin W$,
\begin{equation}\label{modtwo:gap-first-case}
 c_{\mathrm{cond}}-E
 =X-r_0r_1-r+\sum_{i=0}^1
                    \varepsilon_i(n_i-t_i-r_i)
 \geq X-r_0r_1-r.
\end{equation}
The added terms are nonnegative by
\eqref{modtwo:parameter-ranges}. If $e\in W$, then both
$\varepsilon_i=1$, and the exact expression is
\begin{equation}\label{modtwo:gap-second-case}
 c_{\mathrm{cond}}-E=X-r_0r_1-r+1
 \geq X-r_0r_1-r.
\end{equation}
Each summand in $X$ is at least $r_0r_1$, since
$r_i\leq t_i$ and $r_i\leq n_i-t_i$. Consequently
\begin{equation}\label{modtwo:half-cross}
 c_{\mathrm{cond}}-E\geq X/2-r.
\end{equation}

Put $m=\min(t,n-t)$. We have $m\geq1$ and $r\leq m$.
Simultaneously replacing $t_i$ by $n_i-t_i$ does not change $X$.
Thus there are $s_i$, either $s_i=t_i$ for both $i$ or
$s_i=n_i-t_i$ for both $i$, such that $s_0+s_1=m$ and
\[
 X=n_1s_0+n_0s_1-2s_0s_1.
\]
Using balance, $s_0s_1\leq m^2/4$, and $m\leq n/2$, we obtain
\[
 X\geq \frac n3 m-\frac{m^2}{2}\geq \frac{nm}{12}.
\]
Together with \eqref{modtwo:half-cross}, this gives, for $n\geq48$,
\begin{equation}\label{modtwo:uniform-gap}
 c_{\mathrm{cond}}-E
 \geq (n/24-1)m
 \geq nm/48
 \geq n/48.
\end{equation}
The replacement by $s_i$ is only an algebraic estimate for $X$;
it does not assume that the complementary dimensions define
another feasible candidate. In particular, the argument covers
subspaces of dimension greater than $n/2$.

There are at most $(n+1)^4$ choices of $(t_0,t_1,r_0,r_1)$ and
at most four choices of $(\varepsilon_0,\varepsilon_1)$.
Within each feasible nonzero stratum, the union bound,
\eqref{modtwo:two-block-count}, and
\eqref{modtwo:uniform-gap} give
\[
 \Pr[\text{some invariant candidate in the stratum}\mid Ae=d]
 \leq \gamma^4 2^{E-c_{\mathrm{cond}}}
 \leq \gamma^4 2^{-n/48}.
\]
Therefore, uniformly over every balanced $d\in e^\perp$,
\[
 \Pr[M_2(G)\neq\mathbb F_2^n\mid Ae=d]
 \leq 4\gamma^4(n+1)^4 2^{-n/48}
 \leq 1024(n+1)^4 2^{-n/48}.
\]
The bound is uniform in the balanced vector $d$, so averaging over
$d$ and adding \eqref{modtwo:degree-balance} proves
\eqref{modtwo:failure-bound}.
\end{proof}

\begin{remark}\label{modtwo:near-full}
One extreme case illustrates the role of the
$\varepsilon_i$ correction. If both $n_i$ are even and
$W_i=H_i$, then $t=n-2$, $r=2$, and
$\varepsilon_0=\varepsilon_1=1$. This candidate is unique,
both counting exponents $E_i$ vanish, and its conditional
invariance codimension is $n-4$. Omitting the
$\varepsilon_i(n_i-t_i-r_i)$ terms would lose the required
control of precisely such large subspaces.
\end{remark}

\subsection{Comparison with the adjacency walk module}

The prime-two result cannot be deduced from ordinary adjacency
controllability over $\mathbb F_2$, which has a universal parity
obstruction.

\begin{proposition}\label{modtwo:adjacency-half-rank}
For every simple graph on $n\geq1$ vertices,
\[
 \dim_{\mathbb F_2}\mathbb F_2[A]e
 \leq \lceil n/2\rceil.
\]
\end{proposition}

\begin{proof}
Put $\mu_k=e^{\mathsf T}A^ke$. For odd $k=2j+1$,
\[
 \mu_{2j+1}=(A^je)^{\mathsf T}A(A^je)=0
\]
because $A$ is alternating. For even $k=2j$, the identity
$x^{\mathsf T}x=\sum_vx_v$ over $\mathbb F_2$ gives
\[
 \mu_{2j}=(A^je)^{\mathsf T}(A^je)
 =e^{\mathsf T}A^je=\mu_j.
\]
Repeated division of a positive even index by two shows that
$\mu_k=0$ for all $k\geq1$. The positive-degree walk vectors
are therefore pairwise orthogonal and orthogonal to $e$.

If $n$ is even, also $e\cdot e=0$, so the entire walk module
is totally isotropic and has dimension at most $n/2$.
If $n$ is odd, then $e^\perp$ is nondegenerate symplectic
of dimension $n-1$. The positive-degree walk vectors span a
totally isotropic subspace of dimension at most $(n-1)/2$;
adding $e$ proves the bound.
\end{proof}

Thus the degree-parity diagonal removes, with exponentially high
probability, an arithmetic rank obstruction present for adjacency
walks in every graph. The conclusion is cyclicity of $e$; it does
not assert that $\mathbb F_2\langle A,D\rangle$ is the full matrix
algebra.

\subsection{The arithmetic interpretation}

The local criterion of Proposition~\ref{prop:localcriterion} translates
the theorem into an odd-index statement.

\begin{corollary}\label{modtwo:lattice}
For any simple graph $G$, the following conditions are equivalent:
\begin{enumerate}
 \item $M_2(G)=\mathbb F_2^n$;
 \item $\Lambda_G$ has rank $n$ and $I_G$ is odd;
 \item
 $\Lambda_G\otimes_{\mathbb Z}\mathbb Z_{(2)}
 =\mathbb Z_{(2)}^n$, where
 $\mathbb Z_{(2)}=\{a/b\in\mathbb Q:b\text{ is odd}\}$.
\end{enumerate}
Hence each condition holds for $G(n,1/2)$ with the probability
lower bound supplied by Theorem~\ref{modtwo:random}.
On this event $G$ is AD-full over $\mathbb Q$ and $\mathbb R$,
and the finite group $\mathbb Z^n/\Lambda_G$ has odd order.
\end{corollary}

\begin{proof}
Proposition~\ref{prop:localcriterion} gives the equivalence of the
first two conditions. In the same Smith form, odd invariant factors
are precisely the units of $\mathbb Z_{(2)}$, giving the third.
Apply Theorem~\ref{modtwo:random} for the probabilistic statement.
\end{proof}

Recall that a \emph{pointed AD intertwiner} from $G$ to $H$
is an orthogonal matrix $Q$ satisfying
\begin{equation}\label{modtwo:pointed-intertwiner}
 Q\mathbf1=\mathbf1,\qquad
 QA_G=A_HQ,\qquad
 QD_G=D_HQ.
\end{equation}
Its level $\ell(Q)$, when $Q$ is rational, is the least positive
integer for which $\ell(Q)Q$ is integral.

\begin{proposition}\label{modtwo:level}
Suppose $\Lambda_G$ has rank $n$, and let $Q$ be a real
orthogonal matrix satisfying
\eqref{modtwo:pointed-intertwiner} for a simple graph $H$.
Then $Q$ is rational, $\Lambda_H$ has rank $n$, and
\[
 Q\Lambda_G=\Lambda_H,\qquad
 \ell(Q)\mid I_G,\qquad
 \ell(Q)\mid I_H.
\]
In fact $I_G=I_H$.
In particular, if $M_2(G)=\mathbb F_2^n$, every such
intertwiner has odd level.
\end{proposition}

\begin{proof}
The given intertwiner preserves every principal moment.
Since $G$ is AD-full, Theorem~\ref{thm:finiteGNS} identifies $Q$
with the unique rational pointed intertwiner.
Lemma~\ref{lem:arith-pointed} gives full rank and lattice transport,
and Theorem~\ref{thm:arith-level} gives equal indices and the stronger
divisibility $\ell(Q)\mid\gcd(c_G,c_H)$. Each exponent divides
the corresponding index. The last assertion follows from
Corollary~\ref{modtwo:lattice}.
\end{proof}

\begin{corollary}\label{modtwo:random-level}
With probability $1-\exp(-\Omega(n))$, a graph
$G\sim G(n,1/2)$ has no pointed AD intertwiner of even
rational level to any simple graph on $n$ vertices.
The assertion is simultaneous over all competitors and
intertwiners.
\end{corollary}

\begin{proof}
On the event in Theorem~\ref{modtwo:random},
Proposition~\ref{modtwo:level} excludes every even level.
\end{proof}

\begin{remark}[Switching]
\label{modtwo:switching-scope}
The four-vertex Godsil--McKay switching matrix~\cite{GodsilMcKay1982} is
$Q=\operatorname{diag}(\tfrac12J_4-I_4,I_{n-4})$,
up to simultaneous coordinate permutation. It is
orthogonal, fixes $\mathbf1$, and has level two.
The preceding corollary therefore excludes such switches when they
also satisfy $QD_G=D_HQ$; adjacency cospectrality and equality of
degree multisets alone do not imply this identity.

The absence of four-vertex switching sets in a random graph also
follows directly: for a fixed four-set, each outside vertex must
have zero, two, or four neighbours in the set, an event of probability
$((1+6+1)/16)^{n-4}=2^{-(n-4)}$. A union bound gives
$\binom n4 2^{-(n-4)}$. Corollary~\ref{modtwo:random-level}
instead controls every pointed AD intertwiner of even level.
\end{remark}

%% file: article_random_fixed_odd.tex
\section{Fixed odd primes}
\label{oddprime:section}

In this section $p$ is a fixed odd prime. All vector spaces and matrices
are over $\mathbb F_p$. Let $G\sim G(n,1/2)$, let $A$ be its adjacency
matrix, put $e=(1,\ldots,1)^{\mathsf T}$, and set
$D=\operatorname{diag}(Ae)$. We prove a statement stronger than
cyclicity of the principal module.

\begin{theorem}\label{odd:fixed-irreducible}
For every fixed odd prime $p$, there are constants $c_p,K_p>0$ such that
\[
 \Pr\bigl[\text{$A,D$ have a nonzero proper common invariant subspace}\bigr]
 \leq \frac{K_p}{n}+K_p\exp(-c_p n).
\]
Consequently
$\mathbb F_p\langle A,D\rangle e=\mathbb F_p^n$ with probability tending
to one.
\end{theorem}

The constants depend on $p$; uniformity for growing primes is supplied
by Theorem~\ref{growing:uniform}. Symmetry allows us to restrict attention
to invariant subspaces of dimension at most $n/2$. Above a fixed
dimension threshold, the cost of invariance outweighs the number of
subspaces respecting the degree residues. Below that threshold, we
use Wood's code estimate~\cite[Lemma~4.1]{Wood2017}. Its application
requires a separate treatment of subspaces concentrated on a small
set of coordinates.

\subsection{Degree residues and linear constraints}

For $\lambda\in\mathbb F_p$, write
$N_\lambda=|\{i:(Ae)_i=\lambda\}|$.

\begin{lemma}\label{oddprime:balance}
For every fixed $\epsilon>0$,
\[
 \Pr\left[\max_{\lambda\in\mathbb F_p}
       |N_\lambda-n/p|>\epsilon n\right]=O_{p,\epsilon}(n^{-1}).
\]
\end{lemma}

\begin{proof}
Set
\[
 \rho_p=\max_{1\leq j<p}\left|\frac{1+\exp(2\pi i j/p)}2\right|<1.
\]
Fourier inversion shows that the probability of any residue of a
$\operatorname{Bin}(m,1/2)$ variable differs from $1/p$ by at most
$\rho_p^m$. For two distinct vertices, condition on the edge joining
them. Their remaining degree sums are independent
$\operatorname{Bin}(n-2,1/2)$ variables. Thus
\[
 \mathbb E N_\lambda=n/p+O(n\rho_p^{n-1}),\qquad
 \operatorname{Var}(N_\lambda)
       \leq n+O_p(n^2\rho_p^{n-2}).
\]
Chebyshev's inequality and a union bound over the fixed $p$ residues
give the result.
\end{proof}

Call a residue assignment \emph{balanced} if every cell has size in
$[n/p-\epsilon n,n/p+\epsilon n]$. We will sum over deterministic
balanced assignments $g\in\mathbb F_p^n$, while estimating each event
under the original independent-edge distribution.

\begin{lemma}\label{oddprime:constraint}
Let $W\leq\mathbb F_p^n$ have dimension $t$ and radical dimension
$r=\dim(W\cap W^\perp)$ for the standard dot product. Then
$AW\subseteq W$, as a condition on symmetric matrices, has codimension
\[
 c_{\rm sym}=t(n-t)-\frac{r(r-1)}2.
\]
For the adjacency matrix of $G(n,1/2)$,
\[
 \Pr[AW\subseteq W]\leq 2^{-c_{\rm sym}+n}.
\]
\end{lemma}

\begin{proof}
The constraint functionals are the symmetric products of $W$ with
$W^\perp$. Split off their common subspace $W\cap W^\perp$.
The dimension of these products is
\[
 \frac{r(r+1)}2+r(t-r)+r(n-t-r)+(t-r)(n-t-r)
   =t(n-t)-\frac{r(r-1)}2.
\]
Here odd characteristic identifies symmetric tensors with the dual
of symmetric matrices. Restricting to zero-diagonal matrices reduces
the rank of the constraints by at most $n$. A rank-$c$ system of
linear equations in independent Bernoulli variables has probability
at most $2^{-c}$: after fixing the free variables, each of $c$ pivot
variables has a prescribed value, of probability at most $1/2$.
This also proves the displayed bound when its right side exceeds one.
\end{proof}

\subsection{Subspaces above a fixed dimension threshold}

Put
\[
 C_p=\prod_{j=1}^{\infty}(1-p^{-j})^{-1}.
\]

\begin{lemma}\label{oddprime:orthogonal-count}
In a nondegenerate symmetric space of dimension $N$ over $\mathbb F_p$,
the number of $t$-dimensional subspaces of radical dimension $r$ is at most
\[
 C_p^2\,2^r\,p^{\,t(N-t)-r(r+1)/2}.
\]
\end{lemma}

\begin{proof}
First choose a totally isotropic $r$-space $R$. After choosing $j$
independent mutually orthogonal isotropic vectors, the next vector
lies in their orthogonal complement, of dimension $N-j$. The
quadratic form factors through a nondegenerate space of dimension
$N-2j$. At the relevant steps this dimension is at least two, and
the elementary polynomial zero bound gives at most
$2p^{N-j-1}$ choices. Dividing the product of these bounds by the
number of ordered bases of an $r$-space gives
\[
 \#\{R\}\leq C_p\,2^r\,p^{rN-(3r^2+r)/2}.
\]
For fixed $R$, the space $W/R$ has dimension $t-r$ in
$R^\perp/R$, whose dimension is $N-2r$. Discarding its required
nondegeneracy, the Gaussian-binomial bound gives at most
$C_p p^{(t-r)(N-t-r)}$ choices. Multiplication proves the claim.
\end{proof}

Write $\alpha=\log_2p$ and $\theta=1-\alpha/p>0$. Choose $\epsilon>0$
such that
\begin{equation}\label{oddprime:epsilon}
 \epsilon\leq\frac{\theta}{8\alpha},\qquad
 \epsilon<\frac1{4p},\qquad
 \epsilon\log p\leq
 \frac{(1-1/p)\log p-\log2}{4}.
\end{equation}
These conditions are compatible: the last numerator is positive
at $p=3$ and increases with $p$. Fix a deterministic balanced
degree assignment $g$.
Every subspace invariant under $\operatorname{diag}(g)$ decomposes as
$W=\bigoplus_\lambda W_\lambda$. Put
\[
 t_\lambda=\dim W_\lambda,\quad
 r_\lambda=\dim(W_\lambda\cap W_\lambda^\perp),\quad
 t=\sum_\lambda t_\lambda,\quad r=\sum_\lambda r_\lambda,
 \quad B=t(n-t)-r^2/2.
\]
Lemma~\ref{oddprime:orthogonal-count} bounds the number in this stratum
by $C_p^{2p}2^r p^E$, where
\begin{align}
 E&=\sum_\lambda
       \left(t_\lambda(N_\lambda-t_\lambda)
                    -\frac{r_\lambda(r_\lambda+1)}2\right)\notag\\
  &\leq \frac Bp+\epsilon nt-\frac r2.
 \label{oddprime:count-exponent}
\end{align}
The last inequality uses Cauchy--Schwarz separately for the $t_\lambda$
and $r_\lambda$. Combining this with Lemma~\ref{oddprime:constraint}
and at most $p^n$ choices of $g$, the binary logarithm of the union bound,
apart from a constant depending on $p$, is at most
\[
 -\theta B+\alpha\epsilon nt+(\alpha+1)n+\frac{1-\alpha}{2}r.
\]
For $t\leq n/2$ and $r\leq t$, we have $B\geq nt/4$, while
$1-\alpha<0$. Our choice of $\epsilon$ makes the exponent at most
\begin{equation}\label{oddprime:medium-bound}
 -\theta nt/8+(\alpha+1)n.
\end{equation}
Choose a fixed integer
\[
 T\geq 16(\alpha+1)/\theta.
\]
There are at most $(n+1)^{2p}$ dimension-and-radical strata.
Thus the probability, on the balanced event, of a common invariant
subspace with $T\leq t\leq n/2$ is exponentially small in $n$.

\subsection{Codes with a prescribed target}

A matrix $F\in\mathbb F_p^{k\times m}$ is a code of distance $\delta m$
if deleting fewer than $\delta m$ columns leaves a matrix of row rank
$k$. Equivalently, $F$ has row rank $k$ and every nonzero vector in
its row space has support at least $\delta m$.

\begin{lemma}\label{oddprime:code}
Fix an odd prime $p$, an integer $k\geq1$, and $\delta>0$.
Let $F$ be a code of distance $\delta m$ and
let $A_m$ be a zero-diagonal symmetric Bernoulli matrix. For every
fixed target $Y\in\mathbb F_p^{k\times m}$,
\[
 \Pr[FA_m=Y]\leq 2^m K\,p^{-km},
\]
where $K$ depends only on $p,k,\delta$.
\end{lemma}

\begin{proof}
Add an independent Bernoulli diagonal $Z$ to $A_m$. Wood's
arbitrary-target estimate~\cite[Lemma~4.1]{Wood2017} applies to
$X=A_m+Z$, with coefficient modulus $p$, target group
$\mathbb F_p^k$, and balance parameter $1/2$. Indeed, the upper-triangular
entries of $X$, including the diagonal, are independent and have
maximum atom $1/2$ modulo $p$. The code condition is precisely
surjectivity after the allowed column deletions. Hence
$\Pr[FX=Y]\leq Kp^{-km}$, uniformly over $Y$. Independence of $Z$ gives
\[
 \Pr[FA_m=Y]
 =\Pr[FX=Y\mid Z=0]
 \leq 2^m\Pr[FX=Y],
\]
as required.
\end{proof}

We now consider $1\leq t<T$. Fix the code-distance parameter once
for all remaining cases. Put $\beta=1-1/p-\epsilon>0$, and write
$\mathcal H(x)=-x\log x-(1-x)\log(1-x)$ for the binary entropy
with natural logarithms. Choose $q>0$ so that
\begin{equation}\label{oddprime:entropy-choices}
 q<\frac1{4p},\qquad
 2Tq+\frac{\mathcal H(q)}{\log p}<\frac{\beta}{2},\qquad
 \mathcal H(q)+q\log(2p)<\frac{\log2}{2},
\end{equation}
and set $\delta=q/(T+1)$. Such a choice exists because
$\mathcal H(q)\to0$ as $q\to0$.

Fix a deterministic balanced degree assignment $g$ and a $t$-space $W$
invariant under $D_g=\operatorname{diag}(g)$. Encode $W$ in one
of two ways.

If $e\notin W$, let the rows of $U$ be a basis of $W$ and put
$F=[U;e^{\mathsf T}]$, of rank $k=t+1$. There is a matrix $B$ with
$UD_g=BU$, and $B^p=B$. Its eigenspaces $E_\lambda$ form a direct
sum, and the $j$th column $U_j$ belongs to $E_{g_j}$. The event of
interest implies
\[
 UA=C U,\qquad e^{\mathsf T}A=g^{\mathsf T}
\]
for some $t\times t$ matrix $C$. Thus $FA$ has a fixed target once
$F,g,C$ are chosen.

If $e\in W$, take a basis matrix $F$ of $W$ whose last row is
$e^{\mathsf T}$, and put $k=t$. Again $FD_g=BF$ with $B^p=B$,
and $F_j\in E_{g_j}(B)$. Invariance under $A$ gives $FA=CF$.
For an upper bound we may retain its first $k-1$ target rows and
prescribe the last row to be $g^{\mathsf T}$. Every actual instance
satisfies one of these target equations. In both encodings every
column of $F$ has final coordinate one. The numbers of choices of
$B,C$, and of auxiliary subspaces of $\mathbb F_p^k$, are bounded
in terms of $p,T$. We count all such encodings, allowing a subspace
to be represented more than once.

Suppose first that $F$ is a code of distance $\delta n$.
For fixed $g,B$, there are at most
\[
 p^{(k-1)(n/p+\epsilon n)}
\]
possible matrices $F$. In the first encoding this follows from
$\sum_\lambda\dim E_\lambda=k-1$. In the second, the allowed columns
in each $E_\lambda$ lie in the affine slice of final coordinate one.
The directions of these slices are
$E_\lambda\cap\ker(\text{final coordinate})$, whose direct sum has
dimension at most $k-1$.
Lemma~\ref{oddprime:code} and the $p^n$ degree assignments therefore
give exponential factor
\begin{equation}\label{oddprime:code-sum}
 2^n p^{-(k-1)(1-1/p-\epsilon)n}.
\end{equation}
In the first encoding $k\geq2$. In the second, every residue cell
is nonempty on the balanced event; its nonzero columns show that
all $p$ eigenspaces $E_\lambda$ are nonzero. Hence $k\geq p$,
which also excludes $k=1$. Equation~\eqref{oddprime:epsilon} makes
\eqref{oddprime:code-sum} exponentially small.

\subsection{Depth and the zero-diagonal cost}

Suppose $F$ is not a code of distance $\delta n$. Choose the largest
$d\geq1$ such that deleting fewer than $d\delta n$ columns leaves
a span $H\leq\mathbb F_p^k$ of codimension $d$.
Such a $d$ exists because $F$ is not a code. Fix one witness $H$ and let
\[
 h=k-d,\qquad
 \tau=\{j:F_j\notin H\},\qquad
 s=|\tau|,\qquad \eta=[n]\setminus\tau .
\]
The columns indexed by $\eta$ span $H$: they contain every column
left by the witness deletion, and all lie in $H$. Thus
$0<s<d\delta n\leq T\delta n<qn$. Since $q<1/2$, at least one
column remains. Its final coordinate is one, so $h\geq1$.

\begin{lemma}\label{oddprime:depth}
The restricted matrix $F_\eta$, valued in $H$, is a code of distance
$\delta n$. For every fixed target $Y$,
\[
 \Pr[FA=Y]\leq Kp^{-h(n-s)},
\]
where $K$ depends only on $p,T,\delta$.
\end{lemma}

\begin{proof}
If fewer than $\delta n$ additional deletions lowered the rank in
$H$, the total number of deletions would be less than
$(d+1)\delta n$ and their codimension would be some $d'\geq d+1$.
It would also be less than $d'\delta n$, contradicting maximality
of $d$. This proves the code assertion.

Let $\pi:\mathbb F_p^k\to\mathbb F_p^k/H$ be the quotient map.
The target equations on the columns in $\eta$ are
\[
 F_\eta A[\eta,\eta]+F_\tau A[\tau,\eta]=Y_\eta.
\]
Their images under $\pi$ give an event $\mathcal Q$ involving
only the edges between $\tau$ and $\eta$.
At least one coefficient $F_i+H$, $i\in\tau$, is nonzero.
Conditioning all other summands, a Bernoulli sum with such a
coefficient hits a prescribed target with probability at most $1/2$.
Different $j\in\eta$ use disjoint sets of edges, so these quotient
conditions satisfy $\Pr(\mathcal Q)\leq2^{-|\eta|}$.

Let $\mathcal E$ be the information obtained by exposing all edges
meeting $\tau$. Conditional on $\mathcal E$, the matrix $A[\eta,\eta]$
is still a zero-diagonal symmetric Bernoulli matrix. On $\mathcal Q$,
the residual target $Y_\eta-F_\tau A[\tau,\eta]$ lies in $H$.
Choose a basis of $H$ and apply Lemma~\ref{oddprime:code}; the
code distance $\delta n$ is at least $\delta|\eta|$. Its bound is
uniform over the exposed values and the resulting target. Therefore
\[
 \Pr[FA=Y]
 \leq
 \mathbb E\!\left[
   \mathbf1_{\mathcal Q}\,2^{|\eta|}Kp^{-h|\eta|}
 \right]
 \leq Kp^{-h|\eta|}.
\]
The quotient equations exactly pay for the zero-diagonal cost.
\end{proof}

We next count these depth strata while preserving the degree grading.
Fix $B,H,\tau,g$, and put
$H_0=H\cap\ker(\text{final coordinate})$, of dimension $h-1$.
The allowed columns outside $\tau$ in a residue cell $\lambda$ are
\[
 H\cap(E_\lambda\times\{1\})
 \quad\hbox{or}\quad
 H\cap E_\lambda\cap\{\text{final coordinate}=1\},
\]
according to the encoding. Each nonempty set is affine. Its direction
is contained in $H_0$, and the direction spaces for different
$\lambda$ form a direct sum. Their dimensions $a_\lambda$ satisfy
$\sum_\lambda a_\lambda\leq h-1$.
Thus there are at most
\[
 p^{(h-1)(n/p+\epsilon n)+Ts}
\]
choices of $F$: the last term allows arbitrary columns on $\tau$.
The choices of $\tau$ cost $\binom ns$, and the remaining
$B,C,H,d,t$ choices have bounded multiplicity.
After Lemma~\ref{oddprime:depth} and the $p^n$ assignments $g$,
the exponent to base $p$, apart from $\log_p\binom ns$ and a fixed
constant, is at most
\begin{equation}\label{oddprime:depth-count}
 -(h-1)(1-1/p-\epsilon)n+(T+h)s.
\end{equation}
For $h\geq2$, use $h\leq T$, $s<qn$, and
$\sum_{s\leq qn}\binom ns\leq\exp(n\mathcal H(q))$.
Equations~\eqref{oddprime:depth-count} and
\eqref{oddprime:entropy-choices} bound the total, up to a constant
depending on $p,T,\delta$, by
\[
 \exp\!\left(
   n[-\beta\log p+2Tq\log p+\mathcal H(q)]
 \right)
 \leq p^{-\beta n/2}.
\]

\subsection{The one-dimensional residual image}

It remains to treat $h=1$. All columns $F_j$, $j\in\eta$, are equal
to the unique vector $f_0\in H$ of final coordinate one. Since
$s<qn<n/(4p)$ and each residue cell has size greater than
$3n/(4p)$, every cell meets $\eta$.
In the encoding $F=[U;e^{\mathsf T}]$, write $f_0=(u_0,1)$.
The vector $u_0$ then belongs to two distinct eigenspaces of $B$,
so $u_0=0$. Thus $W$ is supported on $\tau$.
In the encoding where $e\in W$, the same argument would force
$f_0=0$, contradicting its last coordinate. That case is impossible.

For fixed $S\subseteq[n]$ of size $s$ and a fixed $t$-space
$W\leq\mathbb F_p^S$, invariance under $A$ forces every row of
$A[[n]\setminus S,S]$ to annihilate $W$. These are independent
rank-$t$ Bernoulli systems. Therefore
\[
 \Pr[AW\subseteq W]\leq 2^{-t(n-s)}.
\]
There are at most $C_p p^{t(s-t)}$ possible $W$ inside $\mathbb F_p^S$.
For each fixed dimension $1\leq t<T$, the total sparse contribution
is therefore bounded by
\[
 C_p\sum_{t\leq s\leq qn}
       \binom ns\,p^{t(s-t)}\,2^{-t(n-s)}.
\]
Discarding the favourable term $-t^2\log p$ and applying the same
entropy bound, its logarithm, apart from $\log C_p$, is at most
\[
 n\bigl[\mathcal H(q)+tq\log p-t(1-q)\log2\bigr]
 \leq -\tfrac12tn\log2,
\]
by \eqref{oddprime:entropy-choices}. The sum over $1\leq t<T$
is exponentially small. This estimate counts all small-support
subspaces directly and incurs no factor for degree assignments.

\begin{proof}[Proof of Theorem~\ref{odd:fixed-irreducible}]
Fix $\epsilon,T,\delta$ as above, in that order.
Degree balance fails with probability $O_p(n^{-1})$.
On balance, \eqref{oddprime:medium-bound} excludes
$T\leq\dim W\leq n/2$, while the code, depth, and sparse estimates
exclude $1\leq\dim W<T$, with exponentially small total probability.
Because $A,D$ are symmetric, the orthogonal complement of every
common invariant subspace is again invariant. Thus every nonzero
proper invariant subspace would yield one of dimension at most $n/2$.
Enlarging $K_p$ covers the finitely many small values of $n$ and proves
the first assertion. The principal module contains
$e\neq0$ and is invariant, so it must be the entire space.
\end{proof}

%% file: article_random_growing_primes.tex
\section{Uniform control of growing primes}
\label{sec:random-growing-primes}

The fixed-prime argument leaves an interval whose upper endpoint grows
with \(n\). We treat it by combining uniform bounds on degree-residue
classes with a dimension split. A small invariant subspace is supported
on few residue classes, which forces a sparse vector in a rectangular
kernel. For larger subspaces, the probability of invariance dominates
the number of candidates. Both estimates use absolute constants.

\begin{theorem}[Simultaneous irreducibility in a growing interval]
\label{growing:uniform}
There is an absolute integer \(P_0\) such that, for
\(G\sim G(n,1/2)\), with probability \(1-\exp(-\Omega(\sqrt n))\) the pair
\((A_G,D_G)\) has no nonzero proper common invariant subspace over
\(\mathbb F_p\), simultaneously for every prime \(P_0\le p\le n\).
Consequently
\[
 M_p(G)=\mathbb F_p^n
 \qquad(P_0\le p\le n)
\]
on that same event.
\end{theorem}

As in the fixed-prime theorem, irreducibility here is over the base
field \(\mathbb F_p\).

\subsection{Uniform sizes of degree-residue classes}

We give the required estimate for integer moduli, although only prime
moduli will be used.  All logarithms with a subscript \(2\) are binary.

\begin{lemma}\label{lem:unimodal-residue-mass}
Let \((a_j)_{j\in\mathbb Z}\) be a finitely supported nonnegative
unimodal sequence of sum \(1\), with maximum \(M\).
For every integer \(p\ge1\) and every residue \(b\pmod p\),
\[
 \sum_{j\equiv b\;(\mathrm{mod}\;p)}a_j\le \frac1p+2M.
\]
\end{lemma}
\begin{proof}
After translating the support, assume it lies in \(\{0,\ldots,N\}\).
Set \(w_j=\mathbf1_{\{j\equiv b\pmod p\}}-1/p\) and
\(W_j=\sum_{i=0}^j w_i\).  The discrepancy in any initial interval
has absolute value at most \(1\), so \(|W_j|\le1\).
Summation by parts, with \(a_{N+1}=0\), gives
\[
 \left|\sum_{j=0}^Nw_ja_j\right|
 \le\sum_{j=0}^N|a_j-a_{j+1}|\le2M.
\]
\end{proof}

\begin{lemma}\label{lem:uniform-degree-residue-size}
For \(n\ge2^{16}\) and \(128\le p\le n\), put
\[
 m_p(G)=\max_{a\in\mathbb Z/p\mathbb Z}
       |\{v:\deg_G(v)\equiv a\pmod p\}|,\qquad
 b_p=\left\lceil32\left(\frac np+\sqrt n\right)\right\rceil .
\]
Then
\[
 \Pr\bigl[\exists\,p\in\mathbb Z,\ 128\le p\le n:
                     m_p(G)\ge b_p\bigr]
 \le n^2(e/8)^{32\sqrt n}.
\]
\end{lemma}
\begin{proof}
Fix \(p\) and a set \(C\) of \(b_p\) vertices.
The stated constants imply \(b_p\le n/2\).
Condition on all edges internal to \(C\).
For \(v\in C\), its external degree is an independent
\(\operatorname{Bin}(n-b_p,1/2)\) random variable; different
vertices of \(C\) use disjoint unordered edges.

The maximum atom of this binomial distribution is at most \(2/\sqrt n\).
Lemma~\ref{lem:unimodal-residue-mass}, applied separately after each
fixed internal-degree shift, gives
\[
 \Pr[\deg_G(v)\equiv a\pmod p\mid G[C]]
 \le\frac1p+\frac4{\sqrt n}.
\]
Thus the probability that every vertex of \(C\) has a specified common
residue is at most \((1/p+4/\sqrt n)^{b_p}\).
Union-bounding over the residue and the set \(C\) yields
\[
 \begin{split}
 \Pr[m_p(G)\ge b_p]
 &\le p\binom n{b_p}
             \left(\frac1p+\frac4{\sqrt n}\right)^{b_p}\\
 &\le p\left[\frac{e(n/p+4\sqrt n)}{b_p}\right]^{b_p}
 \le p(e/8)^{b_p}.
 \end{split}
\]
Sum over all integer moduli \(128\le p\le n\), and use
\(b_p\ge32\sqrt n\).
The bound was obtained for deterministic sets \(C\); averaging over
\(G[C]\) before the union bound preserves its unconditional form.
\end{proof}

\subsection{Counting subspaces and excluding sparse vectors}

\begin{lemma}\label{lem:graded-subspace-count-uniform}
Let a diagonal matrix over \(\mathbb F_p\) have coordinate eigenspaces
of sizes at most \(m\).  For \(1\le t\le n\), the number of its
\(t\)-dimensional invariant subspaces is at most
\[
 (4n)^t p^{mt}.
\]
\end{lemma}
\begin{proof}
Write the coordinate eigenspaces as \(V_1,\ldots,V_q\), with \(q\le n\).
Every invariant subspace decomposes as
\(W=\bigoplus_i W_i\), where \(W_i\le V_i\).
For fixed dimensions \(t_i\) summing to \(t\), the Gaussian-binomial
bound gives
\[
 \prod_{i:t_i>0}
       \genfrac{[}{]}{0pt}{}{\dim V_i}{t_i}_p
 \le 4^t p^{\sum_i t_i(\dim V_i-t_i)}
 \le4^t p^{mt}.
\]
Here
\(\prod_{j\ge1}(1-p^{-j})^{-1}<4\), and at most \(t\) of the
\(t_i\)'s are nonzero.
There are at most \(q^t\le n^t\) possible dimension vectors:
each such vector is realized by at least one ordered allocation of
\(t\) units to the \(q\) classes.
\end{proof}

\begin{lemma}\label{lem:symmetric-invariance-probability}
Let \(p\) be odd, and let \(W\le\mathbb F_p^n\) have dimension
\(1\le t\le n/2\).  For the zero-diagonal Bernoulli adjacency matrix \(A\),
\[
 \Pr[AW\subseteq W]\le 2^{-nt/4+n}.
\]
\end{lemma}
\begin{proof}
Put \(r=\dim(W\cap W^\perp)\).
Lemma~\ref{oddprime:constraint} gives
\[
 \Pr[AW\subseteq W]
 \leq 2^{-c+n},\qquad
 c=t(n-t)-\frac{r(r-1)}2.
\]
This estimate is valid for every odd \(p\), with no constants depending
on \(p\). Since \(r\le t\le n/2\),
\[
 c\ge t(n-t)-\frac{t(t-1)}2\ge\frac{nt}{4}.
\]
\end{proof}

\begin{lemma}\label{lem:uniform-sparse-kernel-vector}
For any prime \(p\) and any integer
\(J\le n/(100\log_2 p)\), the probability that there exist
\(U\subseteq[n]\), \(1\le |U|\le J\), and
\(0\ne x\in\mathbb F_p^U\) such that
\[
 A[[n]\setminus U,U]x=0
\]
is at most
\[
 (n+1)2^{-0.89n}.
\]
\end{lemma}
\begin{proof}
For fixed \(U,x\), every row outside \(U\) contains independent
Bernoulli edge variables.  A nonzero coordinate of \(x\) implies
that its scalar product with such a row is zero with probability
at most \(1/2\).  Different outside rows use disjoint edges.
The event therefore has probability at most \(2^{-(n-|U|)}\).
Union-bounding gives
\[
 \sum_{j=1}^{J}\binom nj p^j2^{-(n-j)}.
\]
Write \(z=j/n\), and let \(H_2(z)\) be binary entropy.
Then \(z\le1/(100\log_2 p)\le0.01\), so
\[
 \begin{split}
 \log_2\left(\binom nj p^j2^{-(n-j)}\right)
 &\le nH_2(z)+j\log_2p-n+j\\
 &\le n\bigl(H_2(0.01)+0.01-1+0.01\bigr)
 <-0.89n.
 \end{split}
\]
Summing over \(j\) proves the assertion.
\end{proof}

\subsection{Proof of the simultaneous theorem}

\begin{proof}[Proof of Theorem~\ref{growing:uniform}]
For \(P_0\le p\le n\), write
\[
 \ell_p=\log_2p,\qquad
 T_p=\lceil16(\ell_p+1)\rceil,\qquad
 b_p=\left\lceil32\left(\frac np+\sqrt n\right)\right\rceil.
\]
Choose \(P_0\ge128\) sufficiently large, and then \(n\) sufficiently
large, that the following inequalities hold uniformly throughout this
interval:
\begin{equation}\label{eq:uniform-prime-parameter-slack}
 T_pb_p\le\frac{n}{100\ell_p},\qquad
 b_p\ell_p\le\frac n{32},\qquad
 \log_2(4n)\le\frac n{32}.
\end{equation}
To justify this order of choices explicitly, use
\[
 \frac{T_pb_p}{n}
 \le561(\ell_p+1)\left(\frac1p+\frac1{\sqrt n}\right).
\]
The supremum of \(\ell_p(\ell_p+1)/p\) for \(p\ge P_0\)
tends to zero as \(P_0\) tends to infinity.
Also
\((\log_2 n)(\log_2 n+1)/\sqrt n\to0\).
These facts establish the first inequality uniformly after the
specified choices. It implies the second since \(T_p\ge1\);
the third holds for all sufficiently large \(n\).

Let
\(\mathcal E=\{m_p(G)<b_p\text{ for every }P_0\le p\le n\}\).
By Lemma~\ref{lem:uniform-degree-residue-size},
\[
 \Pr[\mathcal E^c]\le n^2(e/8)^{32\sqrt n}.
\]

Fix a prime \(p\) in the interval.
Since both \(A\) and \(D\) are symmetric and the ambient dot product
is nondegenerate, the orthogonal complement of a common invariant
subspace is again common invariant.  This remains true when the
subspace itself has a radical.  It therefore suffices to exclude
dimensions \(1\le t\le n/2\).

First consider \(t<T_p\), on the event \(\mathcal E\).
Any \(D\)-invariant \(W\) is graded by the degree-residue classes.
At most \(t\) of those components are nonzero, so \(W\) is supported
on a set \(U\) of size at most
\[
 tb_p\le T_pb_p\le\frac{n}{100\ell_p}.
\]
If \(AW\subseteq W\), any nonzero vector \(x\in W\) satisfies
\(A[[n]\setminus U,U]x=0\).
Lemma~\ref{lem:uniform-sparse-kernel-vector} bounds the probability of
any such obstruction by \((n+1)2^{-0.89n}\).
The sparse-vector estimate already takes the union over all supports
and vectors, so it bounds this case without enumerating degree
assignments.

Now consider \(T_p\le t\le n/2\).
Enumerate the at most \(p^n\) possible degree-residue assignments,
retaining only those with class sizes less than \(b_p\).
For each assignment, Lemma~\ref{lem:graded-subspace-count-uniform}
gives at most \((4n)^t p^{b_pt}\) graded \(t\)-subspaces.
For a fixed such \(W\), the joint event that the assignment is the
actual degree assignment and \(AW\subseteq W\) is contained in the
event \(AW\subseteq W\). Thus the unconditional estimate in
Lemma~\ref{lem:symmetric-invariance-probability} bounds each term,
even though the degree assignment depends on the adjacency matrix.

The binary logarithm of the resulting union bound is at most
\[
 \begin{split}
 n\ell_p+t\log_2(4n)+b_pt\ell_p-\frac{nt}{4}+n
 &\le n(\ell_p+1)-\frac{3nt}{16}\\
 &\le -2n(\ell_p+1),
 \end{split}
\]
using \eqref{eq:uniform-prime-parameter-slack} and \(t\ge T_p\).
Summing over \(t\) costs at most a factor \(n\).

Finally sum over at most \(n\) primes.  The total failure probability
is bounded by
\[
 n^2(e/8)^{32\sqrt n}
 +n(n+1)2^{-0.89n}
 +n^2 2^{-2n(\log_2P_0+1)},
\]
which is \(\exp(-\Omega(\sqrt n))\).
On the complementary event the only common invariant subspaces
are \(0\) and \(\mathbb F_p^n\).
The cyclic module generated by the nonzero vector \(e\)
is common invariant and nonzero, hence equals \(\mathbb F_p^n\).
\end{proof}

\begin{remark}\label{rem:uniform-prime-assembly}
The cutoff \(P_0\) is absolute. The fixed-prime estimates therefore
handle only finitely many remaining primes, while
Theorem~\ref{largeprime:main} already handles all \(p>n\).
\end{remark}

%% file: article_arithmetic_core.tex
\section{Moment fibres and arithmetic intertwiners}
\label{sec:arith-intertwiners}

We now turn from random generation to the arithmetic information
contained in a principal moment function.  The \emph{principal moment
fibre} of \(G\) consists of all graphs \(H\) with \(H\equiv_\mu G\).
On the AD-full locus, Theorem~\ref{thm:finiteGNS} identifies this
relation with pointed rational orthogonal equivalence.  The word
lattice then controls the denominators of the intertwiner.

\begin{lemma}[The unique full pointed isometry]
\label{lem:arith-pointed}
Suppose $G$ is AD-full and $\mu_H=\mu_G$.  Then $H$ is AD-full, and there
is a unique linear isometry
\[
 Q:\mathbb Q^{V(G)}\longrightarrow\mathbb Q^{V(H)}
\]
such that
\[
 Qe_G=e_H,\qquad
 QA_G=A_HQ,\qquad QD_G=D_HQ.
\]
In vertex coordinates, $Q$ is a rational orthogonal matrix, and
$Q\Lambda_G=\Lambda_H$.
\end{lemma}
\begin{proof}
Theorem~\ref{thm:finiteGNS} supplies the unique rational orthogonal
intertwiner.  For every word,
\[
 Qw(A_G,D_G)e_G=w(A_H,D_H)e_H.
\]
These vectors generate the two integer lattices, proving the last
assertion.
\end{proof}

The following is the AD version of the standard walk-matrix denominator
argument.  Bounds using the largest Smith factors on both sides appear,
for example, in \cite{RazaShabbirAbbas2025}.

\begin{theorem}[Index and denominator constraints]
\label{thm:arith-level}
Under the hypotheses of Lemma~\ref{lem:arith-pointed},
\[
                         I_G=I_H.
\]
If $\ell(Q)$ is the least positive integer such that $\ell(Q)Q$ is
integral, then
\[
                  \ell(Q)\mid\gcd(c_G,c_H).
\]
Consequently, if $p\nmid I_G$, then $Q$ is integral over
$\mathbb Z_{(p)}$ and over $\mathbb Z_p$.
\end{theorem}
\begin{proof}
An orthogonal map preserves covolume, so
$Q\Lambda_G=\Lambda_H$ gives equality of their indices in the respective
standard vertex lattices.
By the definition of exponent,
$c_G\mathbb Z^n\subseteq\Lambda_G$.  Hence
\[
 c_GQ\mathbb Z^n\subseteq Q\Lambda_G=\Lambda_H
                         \subseteq\mathbb Z^n.
\]
Thus $\ell(Q)\mid c_G$.  Applying the same argument to $Q^{\mathsf T}$
gives $\ell(Q)\mid c_H$, since transposition does not change the level.
The last assertion follows because every prime dividing $c_G$ divides
$I_G$.
\end{proof}

\begin{remark}
\label{rem:arith-index-not-group}
The index is constant on a full principal moment fibre.
The isometry identifies \(\Lambda_G\) with \(\Lambda_H\), but need
not identify their standard integer ambient lattices.  Thus the
argument gives equality of \(I_G,I_H\), without identifying the
abstract defect groups \(C_G,C_H\).
\end{remark}

\begin{corollary}[One-sided saturation determines the graph]
\label{cor:arith-saturated-rigidity}
If $\Lambda_G=\mathbb Z^{V(G)}$ and $\mu_H=\mu_G$, then $G\cong H$.
No saturation assumption on $H$ and no connectedness assumption is
required.
\end{corollary}
\begin{proof}
Lemma~\ref{lem:arith-pointed} supplies $Q$, and
$Q\mathbb Z^n=\Lambda_H\subseteq\mathbb Z^n$.  Thus $Q$ is an integral
orthogonal matrix, hence a signed permutation matrix.  Since
$Qe=e$, all its signs are positive.
\end{proof}

\begin{corollary}[Random principal-moment rigidity]
\label{integral:random-rigidity}
With probability tending to one, a graph \(G\sim G(n,1/2)\) has the
following property: every finite simple graph \(H\) whose
principal AD moments agree with those of \(G\) through word length
\(2n-1\) is isomorphic to \(G\).
On the same event, every real orthogonal pointed AD intertwiner
from \(G\) to a simple graph is a permutation matrix.
\end{corollary}
\begin{proof}
On the event of Theorem~\ref{integral:random}, \(G\) is integrally
AD-saturated and hence AD-full.  Theorem~\ref{thm:finiteGNS} and
Corollary~\ref{cor:arith-saturated-rigidity} give the first assertion.
The empty-word moment ensures that \(H\) also has \(n\) vertices.
For any real pointed AD intertwiner \(Q\), transport of integer word
combinations gives \(Q\mathbb Z^n=Q\Lambda_G=\Lambda_H
\subseteq\mathbb Z^n\).  Orthogonality and \(Qe=e\) then make \(Q\)
a permutation matrix.
\end{proof}

The conclusion holds simultaneously for all comparison graphs.
Remark~\ref{rem:implications-shorter-rigidity} gives a direct proof
from the integer cut matrices used in the large-prime argument.
The full saturation theorem establishes the stronger integer-generation
property.

\begin{corollary}[Localized switching exclusion]
\label{cor:arith-local-switching}
Suppose $G$ is AD-full and $p\nmid I_G$.  There is no rational orthogonal
matrix of level divisible by $p$ that fixes $e$ and transforms
$(A_G,D_G)$ into the adjacency and degree matrices of a graph.
In particular, saturation at $2$ excludes every such switching matrix
of even level.
\end{corollary}
\begin{proof}
The transformed graph belongs to the same principal moment fibre, so
Theorem~\ref{thm:arith-level} applies.
\end{proof}

For a Godsil--McKay reflection \cite{GodsilMcKay1982} on a switching
set of even size $s\ge4$,
the nontrivial block is $(2/s)J_s-I_s$, whose level is $s/2$.
Consequently, when this switching also intertwines the degree
diagonals, the preceding corollary excludes all sets with $4\mid s$
on a graph saturated at $2$.  Degree preservation is essential here:
the principal AD moments retain the degree operator as well as the
adjacency operator.

\section{Vertex lattices in a principal moment fibre}
\label{sec:arith-fibre}

Fix an AD-full graph $G$ and abbreviate
\[
 V=\mathbb Q^{V(G)},\quad
 A=A_G,\quad D=D_G,\quad \omega=e_G,\quad
 \Lambda=\Lambda_G,\quad I=I_G.
\]
The space $V$ carries its usual positive definite inner product.
All lattices below have full rank in this fixed rational space.
For an integral lattice $K$, write
\[
 K^\#=\{x\in V:\langle x,K\rangle\subseteq\mathbb Z\}.
\]
The finite discriminant group of $\Lambda$ is
\[
 \mathcal D_\Lambda=\Lambda^\#/\Lambda,\qquad
 \beta(x+\Lambda,y+\Lambda)=\langle x,y\rangle+\mathbb Z.
\]
It has order $I^2$ and a nondegenerate symmetric bilinear pairing into
$\mathbb Q/\mathbb Z$.  We use the bilinear pairing, not the discriminant
quadratic form of an even lattice; the vertex lattice need not be even.

Integral unimodular overlattices are described by isotropic subgroups
of the discriminant group.  A walk-lattice application to graph
reconstruction appears in \cite[Section~5]{Xiang2026}.
For the two-operator problem, we must also identify the vertex frame
and require that \(A,D\) become adjacency and degree matrices in that
frame.  These conditions give an exact description of the fibre.

\begin{lemma}[Invariant overlattices and the discriminant module]
\label{lem:arith-overlattices}
The operators $A,D$ preserve $\Lambda^\#$ and induce endomorphisms of
$\mathcal D_\Lambda$.  Integral unimodular lattices $L$ satisfying
\[
             \Lambda\subseteq L\subseteq\Lambda^\#
\]
correspond bijectively to subgroups $S\le\mathcal D_\Lambda$ satisfying
\[
                |S|=I,\qquad \beta(S,S)=0,
\]
by $S=L/\Lambda$.  Moreover, $L$ is stable under $A,D$ if and only if
$S$ is stable under their induced endomorphisms.
\end{lemma}
\begin{proof}
For $x\in\Lambda^\#$ and $y\in\Lambda$,
$\langle Ax,y\rangle=\langle x,Ay\rangle\in\mathbb Z$; the same argument
applies to $D$.  This proves preservation of the dual.
The inverse image of a subgroup $S$ under
$\Lambda^\#\to\mathcal D_\Lambda$ is an intermediate lattice, and all
intermediate lattices arise uniquely in this way.  Its inner products
are integral exactly when $\beta(S,S)=0$.
Since $\operatorname{covol}(\Lambda)=I$, one has
\[
                   \operatorname{covol}(L)=I/|S|.
\]
An integral lattice is unimodular exactly when its covolume is $1$.
This gives the stated cardinality condition.  Stability is equivalent
on the quotient because both operators already preserve $\Lambda$.
\end{proof}

\begin{lemma}[The unit frame is unique]
\label{lem:arith-unit-frame}
Let $L$ be a positive definite integral lattice of rank $n$.
If its vectors of squared norm $1$ span its ambient real space, then
there is an orthonormal basis $v_1,\ldots,v_n$ such that
\[
 L=\bigoplus_{i=1}^n\mathbb Zv_i,\qquad
 \{v\in L:\langle v,v\rangle=1\}=\{\pm v_1,\ldots,\pm v_n\}.
\]
The basis is unique up to signs and ordering.
If a vector $\omega$ has inner product of absolute value $1$ with each
unit vector, exactly one choice of sign on each line satisfies
$\langle\omega,v_i\rangle=1$, and then $\omega=\sum_i v_i$.
\end{lemma}
\begin{proof}
The inner product of two unit vectors is an integer of absolute value
at most $1$.  Equality in Cauchy--Schwarz occurs only when the vectors
are equal up to sign.  Thus distinct unit directions are orthogonal.
The spanning assumption gives exactly $n$ such directions.  Let
$L_0=\bigoplus_i\mathbb Zv_i$.  Integrality gives
$L\subseteq L_0^\#=L_0$, whereas the reverse inclusion is immediate.
The remaining assertions follow by expressing $\omega$ in this
orthonormal basis.
\end{proof}

\begin{definition}[Graph-admissible overlattices]
\label{def:arith-admissible}
An $A,D$-stable integral unimodular overlattice
$\Lambda\subseteq L\subseteq\Lambda^\#$ is \emph{graph-admissible}
if the following conditions hold.
\begin{enumerate}
\item Its unit vectors span $V$.
\item Every unit vector $v\in L$ satisfies
      $|\langle\omega,v\rangle|=1$.
\item On the seed-oriented unit frame
\[
       \mathcal V_\omega(L)
       =\{v\in L:\langle v,v\rangle=1,\ 
                         \langle\omega,v\rangle=1\},
\]
      every vector is an eigenvector of $D$, and
\[
 \langle v,Av\rangle=0,\qquad
 \langle u,Av\rangle\in\{0,1\}\quad(u\ne v).
\]
\end{enumerate}
The graph associated with $L$ has vertex set $\mathcal V_\omega(L)$,
with $u,v$ adjacent when $\langle u,Av\rangle=1$.
\end{definition}

The eigenvector condition on $D$ cannot in general be omitted merely
because $L$ is stable under $D$.  It requires the degree operator to be
diagonal in the prospective vertex frame.  Its diagonal entries then
have the correct degree values automatically: the fixed identity
$A\omega=D\omega$, together with $\omega=\sum_{v\in\mathcal V_\omega(L)}v$,
gives the row sums.

\begin{theorem}[Exact classification of a full principal moment fibre]
\label{thm:arith-fibre}
The isomorphism classes in the principal moment fibre of $G$ are in
bijection with its graph-admissible overlattices
$\Lambda\subseteq L\subseteq\Lambda^\#$.
Equivalently, they are the graph-admissible members of the finite set of
$A,D$-stable isotropic subgroups of $\mathcal D_\Lambda$ of order $I$.
In particular, each such overlattice contributes at most one graph
isomorphism class.
\end{theorem}
\begin{proof}
Let $H$ satisfy $\mu_H=\mu_G$.  Apply
Lemma~\ref{lem:arith-pointed} with the isometry directed from $H$ to $G$,
and denote that map by $R$.  Set
\[
                       L_H=R\mathbb Z^{V(H)}.
\]
It is integral and unimodular, contains
$R\Lambda_H=\Lambda$, and is stable under $A,D$.
Every $x\in L_H$ pairs integrally with $\Lambda\subseteq L_H$, so
$L_H\subseteq\Lambda^\#$.
The images under $R$ of the standard unit vectors form an orthonormal
basis of $L_H$, each with inner product $1$ with $\omega$.  In that basis,
$A,D$ are exactly $A_H,D_H$.  Thus $L_H$ is graph-admissible.

If $H$ is relabelled, composition with the corresponding permutation
gives a pointed intertwiner into $G$.  Uniqueness in
Lemma~\ref{lem:arith-pointed} shows that the resulting lattice is
unchanged.  Hence the construction depends only on the isomorphism
class of $H$.

Conversely, let $L$ be graph-admissible.  By
Lemma~\ref{lem:arith-unit-frame}, its seed-oriented unit frame has
exactly $n$ elements and sums to $\omega$.  Choose an ordering of this
frame and let $R$ have its elements as columns.  Then $R$ is rational
orthogonal and $Re=\omega$.  By admissibility,
$R^{\mathsf T}AR$ is the adjacency matrix of a simple graph $H$ and
$R^{\mathsf T}DR$ is diagonal.  Because $A\omega=D\omega$,
\[
 (R^{\mathsf T}AR)e
       =(R^{\mathsf T}DR)e,
\]
so the latter diagonal is precisely the degree diagonal of $H$.
Every principal moment of $H$ therefore equals that of $G$.
Different orderings of the frame give isomorphic graphs.

The constructions are inverse.  Starting with $H$ recovers its vertex
frame under $R$ and hence its graph up to isomorphism.  Starting with
$L$ gives $R\mathbb Z^n=L$.
The uniqueness of the unit frame is the reason no further multiplicity
of graph classes must be counted within one lattice.
Lemma~\ref{lem:arith-overlattices} gives the subgroup formulation.
\end{proof}

\begin{corollary}[A finite index bound for moment mates]
\label{cor:arith-fibre-bound}
Let $\mathcal F_\mu(G)$ be the set of graph isomorphism classes in the
principal moment fibre of an AD-full graph $G$.  If $I=I_G>1$, then
\[
             |\mathcal F_\mu(G)|
               \le I^{\,2\lfloor\log_2 I\rfloor}.
\]
For $I=1$ the fibre has one member.
\end{corollary}
\begin{proof}
The group $\mathcal D_\Lambda$ has order $I^2$.
Every subgroup of order $I$ can be generated by at most
$k=\lfloor\log_2 I\rfloor$ elements: successively adjoining an element
outside the subgroup already generated at least doubles its order.
Pad a shorter generating tuple with zeros.  The number of such subgroups
is therefore at most the number $(I^2)^k$ of ordered $k$-tuples in
$\mathcal D_\Lambda$.  The graph-admissibility, isotropy, and stability
conditions only decrease the number.  Apply
Theorem~\ref{thm:arith-fibre}.  When $I=1$, the standard vertex lattice
is the sole intermediate unimodular lattice.
\end{proof}

\begin{remark}[Effective scope]
\label{rem:arith-complexity}
The classification is finite and decidable: enumerate the candidate
subgroups and test the corresponding positive definite lattices for
their unit frames and graph admissibility.  The bound counts graph
isomorphism classes; it supplies no running-time bound of the form
\(f(I)n^{O(1)}\) for those tests.
\end{remark}

\section{Small defect and the index-two classification}
\label{sec:arith-nonsaturated}

The fibre description also applies when \(I>1\).  We first record the
elementary restrictions at prime and squarefree index, then classify
the graph realizations completely at index two.

\begin{corollary}[Unique invariant gluing]
\label{cor:arith-unique-gluing}
If $\mathcal D_\Lambda$ has only one $A,D$-stable isotropic subgroup
of order $I$, then $G$ is determined by its principal moments.
More generally, it suffices that only one such subgroup be
graph-admissible.
\end{corollary}
\begin{proof}
The standard vertex lattice is always one graph-admissible member.
Theorem~\ref{thm:arith-fibre} gives the conclusion.
\end{proof}

\begin{proposition}[Prime-index fibres]
\label{prop:arith-prime-index}
Suppose $I_G=p$ is prime.
\begin{enumerate}
\item If $\mathcal D_\Lambda$ is cyclic, then $G$ is determined by its
      principal moments.
\item If $p$ is odd, then $|\mathcal F_\mu(G)|\le2$.
\item If $p=2$, then $|\mathcal F_\mu(G)|\le3$.  If the induced
      nondegenerate bilinear form on
      $\mathcal D_\Lambda\cong\mathbb F_2^2$ is nonalternating, then
      $G$ is determined by its principal moments.
\end{enumerate}
\end{proposition}
\begin{proof}
The discriminant group has order $p^2$.  If it is cyclic, it has a
unique subgroup of order $p$, proving the first assertion.
Otherwise it is elementary abelian of rank two.  Identifying the value
group $p^{-1}\mathbb Z/\mathbb Z$ with $\mathbb F_p$ by
$a/p+\mathbb Z\mapsto a\bmod p$ turns the discriminant pairing into a
nondegenerate symmetric bilinear form over $\mathbb F_p$.
There is at least one isotropic line,
namely $\mathbb Z^{V(G)}/\Lambda$.

For odd $p$, choose a nonzero vector $u$ on that line and a vector $v$
with $\langle u,v\rangle=1$.  Replacing $v$ by
$v-\langle v,v\rangle u/2$ makes $v$ isotropic.
In the resulting basis, the self-pairing of $xu+yv$ is $2xy$.
The only isotropic lines are therefore the two coordinate lines.
There are at most two integral unimodular overlattices, and hence at
most two graph classes.

Over $\mathbb F_2$, there are only three lines in a two-dimensional
space, giving the stated bound.  If the form is nonalternating, the
function $x\mapsto\langle x,x\rangle$ is a nonzero linear functional.
Its kernel is the unique isotropic line.  The first assertion's
uniqueness argument applies again.
\end{proof}

\begin{corollary}[Squarefree index]
\label{cor:arith-squarefree}
Suppose $I_G$ is squarefree and has $s$ odd prime factors.  Then
\[
 |\mathcal F_\mu(G)|
 \le
 \begin{cases}
  2^s,&2\nmid I_G,\\
  3\cdot2^s,&2\mid I_G.
 \end{cases}
\]
\end{corollary}
\begin{proof}
The primary components of the discriminant group are mutually
orthogonal: a pairing between elements of coprime orders is zero.
For each $p\mid I_G$, the $p$-primary component has order $p^2$,
and a unimodular overlattice selects an isotropic subgroup of order
$p$ there.  Apply the local counts in
Proposition~\ref{prop:arith-prime-index} and multiply.
Operator stability and graph admissibility can only reduce the count.
\end{proof}

\begin{theorem}[A binary obstruction that certifies rigidity]
\label{thm:arith-index-two}
Suppose $I_G=2$, and let $v\in\{0,1\}^{V(G)}$ represent the unique
nonzero vector orthogonal modulo $2$ to all AD word vectors.
Then $v$ has even Hamming weight.  If
\[
                         \operatorname{wt}(v)\equiv2\pmod4,
\]
the graph $G$ is determined by its principal moments, although it is
not integrally AD-saturated.
\end{theorem}
\begin{proof}
An index-two subgroup is the kernel of its unique nonzero character
to $\mathbb F_2$, so
\[
             \Lambda=\{z\in\mathbb Z^n:v^{\mathsf T}z
                                            \equiv0\pmod2\}.
\]
Since $\omega=e$ belongs to $\Lambda$, the weight of $v$ is even.
Choose a coordinate $j$ in the support of $v$.
One has
\[
       \Lambda^\#=\mathbb Z^n+\mathbb Z(v/2).
\]
Indeed, the right side pairs integrally with $\Lambda$, contains
$\mathbb Z^n$ with index two, and hence has the same covolume as
$\Lambda^\#$.
The discriminant group is generated by
$a=e_j+\Lambda$ and $b=v/2+\Lambda$.
Both have order two, and they are independent since $v/2$ is not an
integer vector.  In this basis, using the identification
$\tfrac12\mathbb Z/\mathbb Z\to\mathbb F_2$ from the preceding proof,
the bilinear form has matrix
\[
 \begin{pmatrix}
     0&1\\
     1&\operatorname{wt}(v)/2\bmod2
 \end{pmatrix}.
\]
For weight congruent to $2$ modulo $4$, this form is nonalternating.
It has a unique isotropic line, so
Proposition~\ref{prop:arith-prime-index} proves the result.
\end{proof}

\begin{theorem}[Complete classification at index two]
\label{thm:arith-index-two-classification}
Suppose \(I_G=2\), let \(v\) be its unique nonzero binary
obstruction, and put \(S=\operatorname{supp}(v)\).
The principal moment fibre has at most two graph isomorphism classes.
It has exactly two if and only if all of the following hold:
\begin{enumerate}
\item \(|S|=4\);
\item the induced graph \(G[S]\) is regular;
\item every vertex outside \(S\) has zero, two, or four neighbours in \(S\);
\item the four vertices in \(S\) have the same total degree in \(G\).
\end{enumerate}
In that case the unique other member is the Godsil--McKay switch
on \(S\): replace the two neighbours in \(S\) of each outside
vertex having two neighbours by the other two vertices.
If any condition fails, \(G\) is determined by its principal moments.
\end{theorem}
\begin{proof}
Write \(k=|S|\), which is positive and even, and retain the notation
\(a=e_j+\Lambda\), \(b=v/2+\Lambda\) from the preceding proof.
If \(k\equiv2\pmod4\), that proof already gives a unique integral
unimodular overlattice.
If \(4\mid k\), the discriminant form is alternating, and its three
isotropic lines give
\[
 L_a=\mathbb Z^n,\qquad
 L_b=\Lambda+\mathbb Z(v/2),\qquad
 L_{a+b}=\Lambda+\mathbb Z(e_j+v/2).
\]
In the last two lattices, every noninteger vector has half-integer
coordinates at all \(k\) positions of \(S\), and integer coordinates
outside \(S\).  Its squared norm is at least \(k/4\).
When \(k\geq8\), neither lattice has a noninteger unit vector;
its integer unit vectors are only \(\pm e_i\) for \(i\notin S\).
They do not span the space, so both lattices fail graph admissibility.

It remains to consider \(k=4\).  A noninteger unit vector has
coordinates \(\pm\tfrac12\) on \(S\) and zero elsewhere.
In \(L_b\) it has an even number of minus signs, and therefore its
inner product with \(\omega=e\) is \(2,0\), or \(-2\).
Thus \(L_b\) fails the seed-orientation condition.
In \(L_{a+b}\) the number of minus signs is odd.
Its four seed-positive unit vectors on \(S\) have exactly one
minus sign each.  Together with \(e_i\) for \(i\notin S\), they are
the columns of the symmetric orthogonal matrix
\[
 Q_S=\operatorname{diag}\bigl(\tfrac12J_4-I_4,\ I_{n-4}\bigr),
\]
after putting the vertices of \(S\) first.
Consequently this is the only possible other graph realization.

For completeness, graph admissibility here has the asserted
combinatorial form.  Put \(B=A_G[S,S]\) and let \(r_i\) be its
row sums.  The \(i\)th diagonal entry of
\((\tfrac12J_4-I_4)B(\tfrac12J_4-I_4)\) equals
\(|E(G[S])|/2-r_i\).  All these entries vanish exactly when
\(G[S]\) is regular; in that case the conjugation leaves \(B\)
unchanged.  An outside adjacency column \(x\in\{0,1\}^4\)
with coordinate sum \(c\) transforms into \((c/2)e-x\).
This is a zero-one column exactly when \(c\in\{0,2,4\}\).
Finally, if \(d_i\) are the four total degrees, the off-diagonal
entries of the conjugated degree block are
\[
                   \frac14\sum_{s\in S}d_s-\frac12(d_i+d_j).
\]
They vanish for all \(i\ne j\) exactly when all four degrees are
equal.  In that case the degree block is unchanged.  These checks
are necessary and sufficient for the sole candidate lattice to
be graph-admissible, and the outside-column formula is precisely
the stated switch.

When the conditions hold, the resulting graph cannot be isomorphic
to \(G\).  Otherwise a permutation matrix and the nonintegral
matrix \(Q_S\) would both be pointed AD intertwiners from \(G\)
to that graph, contradicting uniqueness in
Lemma~\ref{lem:arith-pointed}.  The exact fibre theorem now
proves all assertions.
\end{proof}

\begin{remark}
Under the index-two hypothesis, the outside-neighbour condition
in the theorem follows automatically from \(|S|=4\).
Indeed the one-dimensional binary obstruction is \(A\)-invariant,
so \(Av\) is a multiple of \(v\); every outside neighbour count
is therefore even.  We retained the condition to make the
connection with the usual switching construction explicit.
The stronger bound of two in this theorem concerns global index
two.  It is not a replacement for the factor three in the
squarefree-index bound with additional odd prime factors.
\end{remark}

\begin{example}[An explicit nonsaturated rigid graph]
\label{ex:arith-index-two}
Let $G$ have vertex set $\{0,\ldots,10\}$ and edge set
\[
\begin{split}
 \{&
 03,06,0\,10,\;
 14,16,17,1\,10,\;
 23,24,25,29,\;
 35,36,37,38,\\
 &45,47,49,4\,10,\;
 57,58,59,\;
 67,6\,10,\;
 78,79,7\,10,\;
 8\,10,9\,10
 \}.
\end{split}
\]
Here $ij$ denotes the unordered pair $\{i,j\}$.
This graph is connected.  Let $v$ be the indicator of
$\{1,2,4,5,7,8\}$.
Direct reduction modulo $2$ gives
\[
               A_Gv=0,\qquad D_Gv=0,\qquad
               v^{\mathsf T}e=0.
\]
Thus every AD word vector lies in the parity kernel of $v$.

Appendix~\ref{app:certificate-index-two} gives three square word
matrices whose determinants have gcd \(2\).  Their nonzero determinants
prove AD-fullness, and their gcd bounds the lattice index by \(2\).
The parity-kernel inclusion forces that index to be even, so \(I_G=2\).
Since $\operatorname{wt}(v)=6$, Theorem~\ref{thm:arith-index-two}
shows that $G$ is determined by all its principal AD moments.
\end{example}

\begin{example}[The index-two bound is attained]
\label{ex:arith-index-two-sharp}
Let \(G\) have vertex set \(\{0,\ldots,9\}\), let
\(S=\{0,1,2,3\}\), and take the edge set
\[
\begin{split}
\{&04,14,25,35,06,36,17,27,08,28,19,39,\\
  &45,49,56,57,58,59,67,68,69,78,89\}.
\end{split}
\]
The set \(S\) is independent, its four vertices have total degree
three, and every outside vertex has two neighbours in \(S\).
Let \(v\) be the indicator of \(S\).  Then \(Av=0\), \(Dv=v\), and
\(v^{\mathsf T}e=0\) modulo two, so
\(\Lambda_G\subseteq\ker(v^{\mathsf T}\bmod2)\).

The two word matrices in Appendix~\ref{app:certificate-sharp} have
determinants with gcd \(2\).  Together with the parity obstruction,
they certify \(I_G=2\).
Theorem~\ref{thm:arith-index-two-classification} now shows that the
fibre has exactly two members.  Its other graph keeps the eleven
edges among \(\{4,\ldots,9\}\) and replaces the neighbours in \(S\)
of each of these six vertices by the complementary pair.
In particular, rational AD-fullness alone does not imply
principal-moment rigidity, and the index-two upper bound is sharp.
\end{example}

\begin{remark}[Two kinds of rigidity]
\label{rem:arith-payoff}
Example~\ref{ex:arith-index-two} separates integral saturation from
principal-moment rigidity, while Example~\ref{ex:arith-index-two-sharp}
shows that rational AD-fullness alone does not force rigidity.
Together they locate the graph-realization question between rational
generation and full integer generation.  Theorem~\ref{integral:random}
shows that a random graph satisfies the strongest of these conditions.
\end{remark}

%% file: article_implications.tex
\section{Further directions}\label{sec:implications}

The ordinary walk lattice provides a useful comparison with the
invariant studied here. By Cayley--Hamilton,
\[
 \Lambda_A(G)=\sum_{j=0}^{n-1}\mathbb Z A^j e
                 \subseteq\Lambda_G.
\]
If \(\Lambda_A(G)\) has full rank, inclusion induces a surjection
\(\mathbb Z^n/\Lambda_A(G)\twoheadrightarrow C_G\). Consequently,
\[
 I_G\mid|\det[e,Ae,\ldots,A^{n-1}e]|,\qquad
 c_G\mid\exp\bigl(\mathbb Z^n/\Lambda_A(G)\bigr).
\]
Closing under \(D\) can therefore remove obstructions retained by
the ordinary walk lattice. The binary theorem makes this distinction
particularly visible: adjacency powers have a modulo-two rank defect
for every \(n>1\), whereas the AD word lattice has none with
exponentially high probability. The all-prime theorem removes the
remaining defect as well.

The proofs suggest the following arithmetic and reconstruction questions.

\paragraph{Other densities and quantitative bounds.}
Determine the range of \(\rho=\rho(n)\) for which
\(\mathbb Z\langle A,D\rangle e=\mathbb Z^n\) with high probability
in \(G(n,\rho)\). The present proof treats \(\rho=1/2\).
For other densities one must control both the degree classes and
integer surjectivity of the relevant cuts; in sparse graphs, repeated
local configurations create additional obstructions.
Even at density \(1/2\), the exponent \(3/20\) reflects the chosen
degree thresholds. Determining the principal causes of failure would
be more informative than optimizing those thresholds alone.

\paragraph{Smaller certificates.}
Proposition~\ref{integral:short-certificate} bounds the necessary
word length by \(n-1\), but the set of all such words is exponential.
The integer construction uses at most two occurrences of \(A\),
at the cost of a common denominator \(\Delta_G^3\).
Can one prescribe, independently of the sampled graph, a
polynomial-size family of words whose evaluations generate
\(\mathbb Z^n\) with high probability? Can integral generation be
achieved with a bounded number of occurrences of \(A\)?
Uniform real spanning families
such as those in~\cite{Lu2026Probes} provide a useful comparison;
integer generation additionally requires the gcd of the maximal
minors to be one.

\paragraph{The least caterpillar diameter.}
Let \(d_*\) be the least fixed integer \(d\) for which the counts
from all caterpillars of diameter at most \(d\) determine
\(G(n,1/2)\) among finite simple graphs with probability tending
to one. Section~\ref{sec:caterpillar-reconstruction} proves
\[
                         5\leq d_*\leq7.
\]
Is \(d_*=5\)? A sufficient strengthening of our construction would be
\[
 \operatorname{span}_{\mathbb Q}
 \{D^ie,D^iAD^je:0\leq i,j<q\}=\mathbb Q^n
 \quad\text{with probability tending to one},
\]
where \(q\) is the number of distinct degrees.
Such a basis, together with the existing rigidity event, would
require moments containing at most three occurrences of \(A\),
and hence test diameter at most five. This spanning condition is
equivalent to full row rank for the neighbour-degree count matrix
\[
 [\,\mathbf1_C,\ A[C,C_1]\mathbf1_{C_1},\ldots,
                    A[C,C_q]\mathbf1_{C_q}\,]
\]
in every degree class \(C\), where \(C_1,\ldots,C_q\) are all
degree classes. The dependence between these classes and the
edge entries is not resolved by the cut estimates proved here.

\paragraph{Realization beyond index two.}
Theorem~\ref{thm:arith-fibre} reduces the principal-moment fibre of
an AD-full graph to a finite realization problem. The index-two
classification solves its first nontrivial case. For larger defects,
one can ask which operator-stable isotropic subgroups admit the
required oriented unit basis, and which graph structures exclude
all but the original vertex lattice. Natural infinite families
with nonzero \(C_G\) and a unique admissible realization would
extend the rigidity example of Section~\ref{sec:arith-nonsaturated}.
A useful next step is a classification at a fixed small odd prime index,
where the unit-vector constraints remain accessible.

For graphs that are not AD-full, the relative torsion group
\((\mathbb Q\Lambda_G\cap\mathbb Z^n)/\Lambda_G\) records only the
arithmetic within the visible module. Additional boundary seeds
may enlarge that module, but a graph reconstruction theorem must
specify how those seeds are chosen and how the resulting local
realizations agree. This is a separate extension of the one-seed
theory developed here.

%% file: article_computation_appendix.tex
\section{Exact computation of the word lattice}
\label{app:integral-computation}

We prove Theorem~\ref{thm:finiteclosure}, retaining the notation of
Section~\ref{sec:finite}.  The argument separates rational stabilization,
index descent, and the size of an exact lattice computation.

\begin{proof}[Proof of Theorem~\ref{thm:finiteclosure}]
The recurrence \(M_{t+1}=M_t+AM_t+DM_t\) implies that an equality
\(M_t=M_{t+1}\) is permanent.  Since \(M_0\) has dimension one,
\(M_{r-1}=M_{\mathbb Q}(G)\).  Choose independent word vectors in
nondecreasing order of length, obtaining an \(n\) by \(r\) matrix
\(B\).  The \(i\)th chosen word has length \(\ell_i\leq i-1\):
otherwise the filtration would have stabilized before acquiring its
required dimension.  The lattice \(\mathbb ZB\) is contained in
\(\Lambda_{r-1}\).  The index of \(\mathbb ZB\) in its saturation is
the gcd of the nonzero \(r\) by \(r\) minors of \(B\), so it bounds
\(I_0\) from above.

Each entry of a length-\(\ell_i\) word vector has absolute value at
most \(\Delta_*^{\ell_i}\).  Hadamard's inequality applied to any
nonzero \(r\) by \(r\) minor gives
\[
 I_0\leq \prod_{i=1}^r
                 \bigl(\sqrt r\,\Delta_*^{\ell_i}\bigr)
       \leq r^{r/2}\Delta_*^{r(r-1)/2}.
\]
For \(t\geq r-1\), all \(\Lambda_t\) have the same rational span.
Every strict enlargement divides their finite saturation index by
an integer at least two.  There can be at most
\(\lfloor\log_2 I_0\rfloor\) such enlargements.
An equality \(\Lambda_t=\Lambda_{t+1}\) makes the lattice stable
under \(A,D\), and hence final.  This proves the stopping bound,
including the case \(I_0=1\).

To compute the lattice, start with the column \(e\) and replace a column
Hermite basis \(B_t\) by a column Hermite basis of
\[
                    [\,B_t\mid AB_t\mid DB_t\,]
\]
until it no longer changes.  This matrix has at most \(3n\) columns.
One additional round after the generation bound certifies stability.

For completeness, the following size bound is independent of the
number of words.  Suppose the word entries at round \(t\) are bounded
by \(B\geq1\) and the current rank is \(s\).  Choose \(s\) independent
word columns \(C\) and \(s\) pivot rows \(P\) on which \(C_P\) is
nonsingular.  Put \(\delta=\det C_P\), so
\(|\delta|\leq s!B^s\).  Projection onto these rows is injective
on the current rational span.  The full-rank projected lattice
contains \(\mathbb ZC_P\), so its index in \(\mathbb Z^s\) is at
most \(|\delta|\).  Its Hermite basis entries have that same bound.
The inverse projection is \(C(C_P)^{-1}\); Cramer's rule bounds
its entries in absolute value by \(s!B^s\).
Lifting the projected Hermite basis gives an integral lattice basis
with entries at most \(s(s!B^s)^2\).  Using the pivot rows of the usual
Hermite procedure gives the same bound for a maintained Hermite basis.
For \(t\leq L+1\), take \(B=\Delta_*^t\).  Since \(L\) in
\eqref{eq:wordlength} is polynomial in \(n\), all these bit lengths are
polynomial in \(n\).  Polynomial-time Hermite and Smith algorithms
\cite{KannanBachem1979} therefore complete the computation in
polynomially many rounds.
\end{proof}

Rational rank can remain constant while the integer lattice strictly
increases.  The index descent is what makes the integral stopping bound
valid; it is separate from the shorter certificate for equality with
\(\mathbb Z^n\) in Proposition~\ref{integral:short-certificate}.

%% file: article_certificates_appendix.tex
\section{Integer certificates for the index-two examples}
\label{app:integer-certificates}

The following determinants make both examples independently
checkable.  In each list, a word is evaluated on the all-ones vector
\(e\), with the rightmost letter acting first.  Rows follow increasing
vertex labels, and columns follow the displayed order.

\subsection{The nonsaturated rigid graph}
\label{app:certificate-index-two}

For the graph in Example~\ref{ex:arith-index-two}, take the ten common
columns
\[
  e,\ Ae,\ A^2e,\ A^3e,\
  DA^2e,\ A^4e,\ ADA^2e,\
  DA^3e,\ A^5e,\ A^2DA^2e.
\]
Appending respectively \(DAe\), \(ADAe\), and \(D^2A^2e\) gives
the three \(11\) by \(11\) matrices with determinants
\[
       -10232927704,\qquad
        14968660956,\qquad
       -778217216902.
\]
Their greatest common divisor is \(2\).  Thus \(\Lambda_G\) has full
rank and its index divides \(2\).  The binary obstruction displayed
in the example proves the reverse divisibility.

\subsection{The graph attaining the index-two bound}
\label{app:certificate-sharp}

For the graph in Example~\ref{ex:arith-index-two-sharp}, use the two
ordered word lists
\[
\begin{split}
\mathcal W_1=(&AAA,\,DAAA,\,AAADA,\,DADAA,\,ADDDAA,\\
             &ADAAADA,\,ADADAAA,\,DAAADAA,\,DDAAAAA,\,DDAAADA),\\
\mathcal W_2=(&AAADA,\,DDDDA,\,AAADAA,\,AADAAA,\,ADDADA,\\
             &AAAADDA,\,AADDDAA,\,ADDDAAA,\,DAADDDA,\,DDAADAA).
\end{split}
\]
The corresponding \(10\) by \(10\) matrices have determinants
\[
\begin{split}
 b_1&=7425047007747097329503430,\\
 b_2&=-339932582914171626171129477243824.
\end{split}
\]
Again \(\gcd(b_1,b_2)=2\).  The word matrices and the binary obstruction
therefore certify \(I_G=2\), as required for the complete classification.